\documentclass{article}

\usepackage{amsmath,amssymb,amsthm}
\usepackage{latexsym}
\usepackage[mathcal]{euscript}
\usepackage{amscd}

\newcommand{\norm}[1]{\| #1 \|}
\newcommand{\origin}{o}
\newcommand{\Harmonics}{\mathcal{H}}
\newcommand{\C}{\mathbf C}

\newcommand{\opi}{{}^{\scriptscriptstyle 0}\pi}
\DeclareMathOperator{\im}{im}
\DeclareMathOperator{\Real}{Re}
\DeclareMathOperator{\Dom}{Dom}

\newtheorem{Theorem}{Theorem}[section]
\newtheorem{corollary}[Theorem]{Corollary}
\newtheorem{lemma}[Theorem]{Lemma}
\newtheorem{prop}[Theorem]{Proposition}

\numberwithin{equation}{section}

\title{Deformation theory of five-dimensional CR %
       structures and the Rumin complex%
\footnote{2000 Mathematics Subject Classification: Primary 32G07, Secondary 32S30, 32V20 }}
\author{Takao Akahori ({\sl akahorit@sanynet.ne.jp}),\\ 
Peter M. Garfield ({\sl garfield@math.washington.edu}),\\ 
and John M. Lee ({\sl lee@math.washington.edu})}

\begin{document}

\maketitle
\begin{abstract}
We construct a versal family of deformations of CR
structures in five dimensions, using a differential complex
closely related to the differential form complex 
introduced by Rumin for contact manifolds. 
\end{abstract}

\section{Introduction}

A natural problem in several complex variables is that of
classifying the deformations 
of an isolated singularity in a complex-analytic variety.
The problem is solved by constructing a 
``versal family'' of deformations of 
the singularity, which is, roughly speaking, 
a minimal family of deformations that includes biholomorphic
representatives of all other
deformations.  (See Section \ref{sect:versality}
for a precise definition.)

Versal families for isolated singularities
were first constructed from an algebraic point of 
view in the late 1960s and early 1970s by Tjurina,
Grauert, and Donin 
\cite{Tjurina:1969,Grauert:1972,Donin:1972}.
Shortly thereafter, 
M. Kuranishi \cite{Kuranishi:1977} outlined a 
program for relating deformations of an isolated
singularity to deformations of the CR structure on a
real hypersurface obtained
by intersecting the variety with a small sphere surrounding
the singular point (the ``link'' of the singularity).  
Then 
Kuranishi's construction was extended and simplified
by subsequent work of the first author and others
\cite{Akahori:1981,Akahori:1982,Akahori:1989,Miyajima:1990,Miyajima:1991,Buchweitz-Millson:book}.

A fundamental limitation of all of these results has
been a dimensional restriction: Because the
deformation complex that was introduced in 
\cite{Akahori:1981,Akahori:1982,Akahori:1989}
failed to be subelliptic in low dimensions,
these results only applied to CR manifolds 
of dimension $7$ or more (and therefore
to singularities of varieties whose 
complex dimension is at
least $4$).

The purpose of this paper is to 
extend the Kuranishi
construction
of versal families of CR structures
to the case of $5$-dimensional CR manifolds.
The new idea here is
a subelliptic estimate and consequent
Hodge theory for a certain subcomplex
of the standard deformation complex inspired by recent
work of M. Rumin on contact manifolds.

Recently, Miyajima \cite{Miyajima:1999}
introduced an alternative approach to constructing
versal families in all dimensions, 
based on analyzing deformations
not only of the CR structure, but of the CR structure
together with its embedding into $\C^N$.
The present approach is of independent interest,
however, because it represents a completion of the
original Kuranishi program of 
constructing an intrinsically-defined versal
family of deformations of the CR structure itself.  
There appears to be little hope for extending
this intrinsic approach to the case of 
$3$-dimensional CR manifolds,
because the relevant cohomology groups in that case
are infinite-dimensional.

Let $(M,{}^0T^{\prime\prime})$ be a compact strictly pseudoconvex CR
manifold of real dimension $5$.  
Deformations of the CR structure of $M$ can be represented as 
$T'$-valued $(0,1)$-forms, where $T'$ is a $3$-dimensional
complex subbundle of $\C\otimes TM$ transverse to the antiholomorphic
tangent bundle ${}^0T''$ (see Section 2 for precise definitions).  
The space of such forms fits into a complex
$(\Gamma(M,T'\otimes \wedge^j(^0T^{\prime\prime})^{\ast},
\overline\partial_{T'}^{(j)})$, the
{\it standard deformation complex}
\cite{Akahori:1981,Buchweitz-Millson:book}.  
In earlier work on higher-dimensional CR deformation theory,
the first author 
defined a subcomplex
$(\Gamma(M,E_j),\overline \partial _j)$ of the standard deformation complex
corresponding to deformations of the CR structure that leave
the contact structure fixed.  
When 
$\dim M = 2n-1 \geq 7$, there is a subelliptic estimate on 
$\Gamma(M,E_2)$, which leads to the construction
of a versal family
\cite{Akahori:1981,Akahori:1982}.  
But if $\dim M = 5$, there is no such estimate.

In this paper, inspired by the differential-form complex introduced by Rumin
\cite{Rumin:1994} for studying de Rham theory on contact manifolds,
we extend the $E_i$ complex by defining a new second-order operator $D$:
$$
0 \rightarrow 
	\Gamma(M,F) \stackrel{D}\rightarrow 
	\Gamma(M,E_1) \stackrel{\overline\partial_1}\rightarrow 
	\Gamma(M,E_2),
$$
where $F$ is a one-dimensional subbundle of $\C \otimes TM$ transverse
to ${}^0T^{\prime\prime} \oplus \overline{{}^0T^{\prime\prime}}$.  
This is closely related to Rumin's complex, in a way we will explain in
Section 4.  A similar complex has also been used in 
\cite{Buchweitz-Millson:book}.  

Once we have proved an 
a priori estimate on $\Gamma(M,E_1)$,
it follows that there is a
Kodaira-Hodge
decomposition theorem on $\Gamma(M,E_1)$. Using
techniques similar to those in 
\cite{Akahori:1981,Akahori:1982}, this leads to a
construction of the versal family in the $5$-dimensional case. 
We remark that Rumin has recently suggested a simpler proof
of an analogous estimate for the complex version of his
complex in arbitrary dimensions.  We hope to pursue this
further in another paper.

\section{Background and Notation}
Let $(M,{}^0T^{\prime\prime})$ be a CR manifold.  By this we mean that
$M$ is a smooth manifold of dimension $2n-1$ and $^0T^{\prime\prime}$ is a
complex subbundle of the complexified tangent bundle $\C \otimes TM$
satisfying  
\begin{align*}
  & {}^0T^{\prime\prime}\cap \overline{{}^0T^{\prime\prime}}=0, 
    \ \mbox{dim}_{\C} {}^0T^{\prime\prime} = n-1, \\
  & [X,Y] \in \Gamma(M,{}^0T^{\prime\prime})
    \mbox{ for all $X,Y \in \Gamma(M,{}^0T^{\prime\prime})$},
\end{align*}
where by $\Gamma(M,E)$ we mean the space of
$C^{\infty}$ sections of the bundle $E$.  
For convenience we will write $^0T^{\prime}$ for 
$\overline{^0T^{\prime\prime}}$ and $H$ for the real bundle 
$\mbox{Re}(^0T^{\prime\prime} \oplus {}^0T^{\prime})$.  
We assume that there is a global non-vanishing real one-form $\theta$
that annihilates $H$; that is, such that $\theta(X) =
\theta(\overline{X}) = 0$ for all $X \in
\Gamma(M,{}^0T^{\prime\prime})$.   Since $H$ is naturally oriented,
the existence of such a form is equivalent to $M$ being orientable.

We define the \emph{Levi form} $L_{\theta}$ by 
\begin{equation}
L_{\theta}(X,\overline Y)
 = - i \theta([X,\overline Y]) 
   \ \text{for} \ X,Y \in {}^0T^{\prime}.  
\end{equation}
If this Levi form $L_{\theta}$ is positive definite or negative
definite, then $(M,{}^0T^{\prime\prime})$ is called 
\emph{strictly} (or \emph{strongly}) \emph{pseudoconvex}.  (After this
section, we will always assume that our CR structure is strictly
pseudoconvex.)  Notice that the Levi form gives us a metric
on $H = \Real ({}^0T^{\prime} \oplus {}^0T^{\prime\prime})$ that
extends to a Riemannian metric on all of $TM$ by declaring that $\xi$
is unit length and orthogonal to $H$.  We will call this metric the
\emph{Webster metric} (see \cite{Webster:1978}). 

When $(M,{}^0T^{\prime\prime})$ is 
strictly pseudoconvex, 
we will call a choice of $1$-form
$\theta$ a \emph{pseudohermitian structure}.  
Let $\xi$ be the unique real vector field satisfying 
$\theta(\xi) = 1$ and $d\theta(\xi,X) = 0$ for all $X\in H$.
Notice that this implies that for every point $p$ of $M$,
$\xi _p \notin \C \otimes H_p = {}^0T_p^{\prime\prime}+{}^0T_p^{\prime}$.

Let $F$ denote the complex line bundle $\C \xi$, and set
$T' := {}^0T^{\prime} + \C\xi$.  We then get vector bundle
decompositions 
\begin{align}
\label{eqn:vector-bundles-1}
\C TM & = T'+{}^0T^{\prime\prime} \\
\intertext{and}
\label{eqn:vector-bundles-2}
\C TM & = {}^0T^{\prime}+{}^0T^{\prime\prime}+F.
\end{align}
Note that these decompositions depend on the choice of $\theta$ (and
thus $\xi$) and so are not CR-invariant.   We will often take
advantage of these decompositions to project onto various components.
For a vector $X$, let us write $\pi_F(X)$ for the $F$-component of
$X$, $\pi'(X)$ for the $T'$-component, $\opi'(X)$ for the
${}^0T'$-component, and $\opi''(X)$ for the ${}^0T''$-component,
according to these decompositions.  Moreover, since we will often be
dealing with vector-valued forms, let us use the same notation for the
projection of, say, 
$\C \otimes TM \otimes \Lambda^j ({}^0T^{\prime\prime})^*$
into component parts  
$F \otimes \Lambda^j ({}^0T^{\prime\prime})^*$,
$T^{\prime} \otimes \Lambda^j ({}^0T^{\prime\prime})^*$,
${}^0T^{\prime} \otimes \Lambda^j ({}^0T^{\prime\prime})^*$,
and
${}^0T^{\prime\prime} \otimes \Lambda^j ({}^0T^{\prime\prime})^*$
via equations \eqref{eqn:vector-bundles-1} and \eqref{eqn:vector-bundles-2}.

It is often useful to identify $\C \otimes \Lambda^k M$ with 
$\C \otimes \Lambda^k H^* \oplus \theta \wedge \C \otimes \Lambda^{k-1}
H^*$.  Notice that this identification depends on the choice of
$\theta$, and so is not CR-invariant.  There is a natural bigrading on
$\C \otimes \Lambda^k H^*$, so we may make a further identification
\begin{equation}
\label{eqn:bigrading}
\C \otimes \Lambda^k M
 = \sum_{p+q=k} \Lambda^{p,q}H^*
   + \theta \wedge \sum_{p+q=k-1} \Lambda^{p,q}H^*.
\end{equation}
This allows us to identify, for example,
$\Lambda^q ({}^0T^{\prime\prime})^* = \Lambda^{0,q}H^*$ 
with honest forms on $M$.

Finally, we note that we will use the Einstein summation convention
whenever possible.  We will use Roman indices ($j$, $k$, for example)
to indicate sums from $1$ to $2n-1$, and Greek indices ($\alpha$,
$\beta$, and so on) for sums from $1$ to $n-1$.

\section{Review of {C}{R} deformation theory}
In this section we survey previous work on the deformation theory of
CR structures.  This work was initiated by Kuranishi
\cite{Kuranishi:1977} as a CR analogue of his work on complex
manifolds.  Most of the work reviewed here was done by the first
author 
\cite{Akahori:1978-delbarb,Akahori:1978-family,Akahori:1981,Akahori:1982}

Following work of the first author \cite{Akahori:1978-family}, we
introduce a first order differential operator 
$\overline\partial_{T'} : \Gamma(M,T') \to \Gamma(M,T'\otimes
(^0T^{\prime\prime})^{\ast})$ 
by 
\begin{equation}
\overline\partial_{T'} Y(\overline X)= \pi' [\overline X,Y]
\qquad \mbox{ for $Y \in \Gamma(M,T')$ and
		 ${\overline X}\in \Gamma(M,{}^0T^{\prime\prime}$).}
\end{equation}
As in the case of scalar-valued differential forms, this generalizes
to operators 
$\overline\partial^{(p)} :
 \Gamma(M,T' \otimes \Lambda^p (^0T^{\prime\prime})^*) 
 \to \Gamma(M,T' \otimes \Lambda^{p+1} (^0T^{\prime\prime})^*) 
$
($p=1,2,\cdots$) given by
\begin{eqnarray}
\lefteqn{
\overline\partial^{(p)} \phi (\overline X_{1}, \ldots, \overline
X_{p+1})
} \nonumber \\
   & = &
	\sum_{j=1}^{p+1} 
			(-1)^{j+1} \pi' [ \overline X_{j},
			\phi(\overline X_1, \ldots, \widehat{\overline
			X}_{j}, \ldots,
				 \overline X_{p+1})] 
	\label{eqn:delbarphi-defn} 
	\\
     & & 
	+ \sum_{j<k} (-1)^{j+k}
			\phi( [\overline X_j, \overline X_k],
			\overline X_1, \ldots,
			\widehat{\overline X}_{j}, \ldots,
			\widehat{\overline X}_{k},
			\ldots, \overline X_{p+1} ) 
	\nonumber 
\end{eqnarray}
for $\phi \in \Gamma(M,T' \otimes \Lambda^p (^0T^{\prime\prime})^*)$
and $\overline X_k \in \Gamma(M,{}^0T^{\prime\prime})$.  We then have
a differential complex  
\begin{equation}
\label{eqn:standard-complex}
\begin{aligned}
  0  \to  \Gamma(M,T') 
     \stackrel{\overline\partial_{T'}}{\to} 
	\Gamma(M,T'\otimes (^0T^{\prime\prime})^{\ast}) 
     \stackrel{\overline\partial^{(1)}}{\to} 
        \Gamma(M,T'\otimes\Lambda^2 (^0T^{\prime\prime})^{\ast}) 
     \stackrel{\overline\partial^{(2)}}{\to} \\
     \cdots \to \Gamma(M,T'\otimes\Lambda^p (^0T^{\prime\prime})^{\ast}) 
     \stackrel{\overline\partial^{(p)}}{\to}
         \Gamma(M,T'\otimes\Lambda^{p+1}(^0T^{\prime\prime})^{\ast}) 
     \to \cdots
\end{aligned}
\end{equation}
with $\overline\partial^{(p+1)}\overline\partial^{(p)}=0$  (see
\cite{Akahori:1978-family}).  
This complex is called the \emph{standard deformation complex}.

A complex subbundle $E \subset \C \otimes TM$ is an \emph{almost CR 
structure} (and the pair $(M,E)$ is an \emph{almost CR manifold}) if
$E\cap \overline E = 0$ and $\mbox{dim}_{\C} E = n-1$.  An almost CR
structure $E$ is \emph{at finite distance from ${}^0T^{\prime\prime}$}
if $\opi^{\prime\prime}|_E : E \to {}^0T^{\prime\prime}$ is a bundle 
isomorphism.  These almost CR structures are characterized by the fact 
that they are graphs over ${}^0T^{\prime\prime}$: there is a bijective
correspondence between elements 
$\phi \in \Gamma(M,\mathrm{Hom}({}^0T^{\prime\prime},T'))
 = \Gamma(M,T' \otimes ({}^0T^{\prime\prime})^{\ast})$
and almost CR structures
\begin{equation*}
{}^{\phi}T^{\prime\prime} 
 := \left\{ \overline X + \phi(\overline X) \ : 
    \ \overline X \in {}^0T^{\prime\prime} \right\}
\end{equation*}
at finite distance from ${}^0T^{\prime\prime}$ (see, for example,
\cite[Proposition 1.1, page 618]{Akahori:1978-delbarb}).  The almost
CR structure ${}^{\phi}T^{\prime\prime}$ is a CR structure exactly
when it satisfies the integrability condition, which can be written as
the non-linear partial differential equation  
\begin{equation*}
P(\phi) := \overline\partial^{(1)}\phi + R_2(\phi) + R_3(\phi) 
	= 0,
\end{equation*}
where $R_k(\phi) \in \Gamma(M,T^{\prime} \otimes \Lambda^2
({}^0T^{\prime\prime})^{\ast})$ ($k=2,3$) are the parts of $P(\phi)$
that are degree $k$ in $\phi$.  They are given by 
\begin{equation}
\label{eqn:R2-defn}
R_2(\phi)(\overline X,\overline Y)
 = \pi' [\phi(\overline X),\phi(\overline Y)]
   - \phi( \opi^{\prime\prime} [\overline X,\phi(\overline Y)] +
           \opi^{\prime\prime} [\phi(\overline X),\overline Y] )
\end{equation}
and
\begin{equation}
\label{eqn:R3-defn}
R_3(\phi)(\overline X,\overline Y)
 = - \phi( \opi^{\prime\prime} [\phi(\overline X),\phi(\overline Y)]).
\end{equation}
See \cite[Theorem 2.1, page 619]{Akahori:1978-delbarb} and the proof
given therein for details.

If we consider only deformations $\phi$ that preserve the contact
structure (that is, for which
${}^{\phi}T^{\prime\prime} \oplus \overline{{}^{\phi}T^{\prime\prime}}
 = {}^{0}T^{\prime\prime} \oplus \overline{{}^{0}T^{\prime\prime}}$),
then we are simply restricting to 
$\phi \in 
\Gamma(M,{}^0T^{\prime} \otimes ({}^0T^{\prime\prime})^*)$.
For such $\phi$, we notice that $R_3(\phi) = 0$ and that 
$\pi_F R_2(\phi) = 0$ (so $\opi' R_2(\phi) = R_2(\phi)$).  Thus 
$P(\phi) = \pi_F \overline \partial^{(1)} \phi 
	+ \opi' \overline \partial^{(1)} \phi 
	+ R_2(\phi)$.  
Our integrability condition $P(\phi) = 0$ is thus equivalent in this
case to $\pi_F \overline \partial^{(1)} \phi = 0$ and 
$\opi' \overline \partial^{(1)} \phi + R_2(\phi) = 0$.  
(Compare \cite[Proposition 1.7.3, page 797]{Akahori:1978-family}.)
This in part motivates the definition of the following subspaces of
$\Gamma(M,{}^0T' \otimes
		    \Lambda^p(^0T^{\prime\prime})^{\ast})$:
\begin{equation}
\Gamma_p = \{ u \in \Gamma(M,{}^0T' \otimes
		    \Lambda^p(^0T^{\prime\prime})^{\ast})
		: \pi_F \overline \partial^{(p)} u = 0 \}. 
\end{equation}
For $\phi \in \Gamma_1 \subset 
	\Gamma(M,{}^0T' \otimes (^0T^{\prime\prime})^{\ast})$,
then, the integrability condition becomes 
$P(\phi) = \opi' \overline \partial^{(1)} \phi + R_2(\phi) = 0$.

We remark that contrary to appearances, 
the definition of $\Gamma_p$ is an algebraic condition
on $u$, not a differential one.  
To see this, apply the one-form
$\theta$ to both sides of equation~\eqref{eqn:delbarphi-defn}.  By the
definition of $\Gamma_p$, the left-hand side is zero, and so 
\begin{align*}
0 
   = & 
	\sum_{j=1}^{p+1} 
			(-1)^{p+1} \theta( [ \overline X_{j},
			u(\overline X_1, \ldots, \widehat{\overline
			X}_{j}, \ldots,
				 \overline X_{p+1})] 
				)  \\
     & 
	+ \sum_{j<k} (-1)^{j+k}
			\theta( 
			u( [\overline X_j, \overline X_k],
			\overline X_1, \ldots,
			\widehat{\overline X}_{j}, \ldots,
			\widehat{\overline X}_{k},
			\ldots, \overline X_{p+1} ) 
				)
\end{align*}
Since $u$ maps into ${}^0 T^{\prime}$, 
which is annihilated by $\theta$,
the second sum is a sum of zeros.  Using
$\theta([X,Y]) = - d\theta(X,Y)$ for $X,Y \in \C \otimes H = 
{}^0 T^{\prime} \oplus {}^0 T^{\prime\prime}$, the first sum becomes  
\begin{equation}
\label{eqn:Gammaj-alg-cond}
0 = \sum_{j=1}^{p+1} (-1)^p
		d\theta( 
			\overline X_{j},
			u(\overline X_1, \ldots, \widehat{\overline
			X}_{j}, \ldots, \overline X_{p+1})
			).
\end{equation}
This is an algebraic condition on $u$.

In fact, the spaces $\Gamma_p$ are smooth sections of
vector bundles.  There are 
\cite[Proposition 2.1, page 313]{Akahori:1981} 
subbundles 
$E_p \subset T^{\prime} \otimes \Lambda^p({}^{0}T^{\prime\prime})^{\ast}$
such that $\Gamma_p = \Gamma(M,E_p)$.  By restricting
$\overline{\partial}^{(p)}$ to $E_p$, we get a sequence of maps
$\overline{\partial}_p$ 
\begin{equation}
\label{eqn:Ej-complex-1}
0 \to \Gamma(M,E_0)
  \stackrel{\overline{\partial}_0}{\longrightarrow} \Gamma(M,E_1)
  \stackrel{\overline{\partial}_1}{\longrightarrow} \Gamma(M,E_2)
  \stackrel{\overline{\partial}_2}{\longrightarrow} \Gamma(M,E_3)
  \stackrel{\overline{\partial}_3}{\longrightarrow} \cdots
\end{equation}
and ${}^{\phi}T^{\prime\prime}$ is integrable for $\phi \in \Gamma_1$
if and only if $P(\phi) = \overline \partial_1 \phi + R_2(\phi) = 0$.

It turns out that $E_0 = 0$ and the resulting complex 
\begin{equation}
\label{eqn:Ej-complex-2}
0 \longrightarrow 0
  \longrightarrow \Gamma(M,E_1)
  \stackrel{\overline{\partial}_1}{\longrightarrow} \Gamma(M,E_2)
  \stackrel{\overline{\partial}_2}{\longrightarrow} \Gamma(M,E_3)
  \stackrel{\overline{\partial}_3}{\longrightarrow} \cdots
\end{equation}
is a differential subcomplex of the standard deformation complex
 (see \cite[Theorem 2.2, page 314]{Akahori:1981}).  
This subcomplex still contains enough information
to be useful; for example, the
inclusion map $\iota : \Gamma(M,E_p) \to 
\Gamma(M,{}^{0}T^{\prime} \otimes \Lambda^p({}^{0}T^{\prime\prime})^{\ast})$
induces a map
\begin{equation*}
\iota^{\ast} : 
         \frac {\ker \overline\partial_p}
	       {\im \overline\partial_{p-1}}
	 \longrightarrow
         \frac {\ker \overline\partial^{(p)}}
	       {\im \overline\partial^{(p-1)}}
\end{equation*}
that is an isomorphism if $p \ge 2$ and surjective if $p=1$
\cite[Theorem 2.4, page 315]{Akahori:1981}.

Furthermore, there are a subelliptic estimate for this complex
\cite[Theorem 4.1, page 319]{Akahori:1981} and a Kodaira-Hodge
decomposition theorem  for $\Gamma(M,E_2)$ \cite[Theorem 4.1, page
328]{Akahori:1981}, provided $\dim M = 2n-1 \ge
7$.  That is, if we define the Laplacian 
$\square = \overline \partial_2^* \overline \partial_2 
        + \overline \partial_1 \overline \partial_1^*$,
then there is a harmonic projector $H$ such that $\square H u = 0$ for
all $u \in \Gamma(M,E_2)$ and a Neumann operator $N$ such that $NH u =
HN u = 0$ and $u = \square N u + H u$ for all $u \in \Gamma(M,E_2)$.  
This construction fails if $\dim M = 5$, as there is no
subelliptic estimate for this complex.

\section{The new complex}
In this section, we introduce a new complex as a replacement for the
differential subcomplex \eqref{eqn:Ej-complex-2} of the standard
differential complex.  
Set 
\begin{equation}
H_0 = \{v\in \Gamma(M,T')
	: \pi_F \overline\partial_{T'} v = 0 \}. 
\end{equation}
We then get a new differential subcomplex of the standard
differential complex \eqref{eqn:standard-complex}:
\begin{align}
\label{eqn:H0-complex}
0 \longrightarrow H_0
  \stackrel{\overline{\partial}_0}{\longrightarrow} \Gamma(M,E_1)
  \stackrel{\overline{\partial}_1}{\longrightarrow} \Gamma(M,E_2)
  \stackrel{\overline{\partial}_2}{\longrightarrow} \cdots.
\end{align}
This complex is a generalization of ideas of the first author that is
new for use in this setting, but it has been introduced by  
Buchweitz and Millson \cite[page 82]{Buchweitz-Millson:book} based in part
on ideas of the third author.
It is straightforward to see that this is a complex: 
the definition of $H_0$ ensures that 
$\overline\partial_0 u \in 
 \Gamma(M,{}^0T'\otimes(^0T^{\prime\prime})^{\ast})$
and the fact that \eqref{eqn:H0-complex} is a subcomplex of the
standard differential complex \eqref{eqn:standard-complex} means that,
in fact, $\overline\partial_0 u \in \Gamma(M,E_1)$.  

We would like to make a few remarks about $H_0$.  
It is not the space of smooth sections of a vector bundle over $M$;
rather, it is the image of a first order differential operator.  We
define this operator $\Gamma(M,F) \to \Gamma(M,T')$ as follows: 
for $Z \in \Gamma(M,F)$, we may write $Z = u\cdot \xi$ for some smooth
function $u$ (namely, $u = \theta(Z)$).  We then get an element 
$X_u \in \Gamma(M,{}^0T')$ by requiring that $u\xi + X_u \in H_0$:
$\pi_F \overline\partial_{T'} (X_u + u\xi) = 0$.  
This is equivalent to $\theta([\overline Y,X_u+u\xi]) = 0$ for
all $\overline Y \in \Gamma(M,{}^0T^{\prime\prime})$.   Another way to
write this is 
\begin{equation}
\label{eqn:Xu-equation}
d\theta(\overline Y,X_u) = {\overline Y} u,
\end{equation}
because $\theta(\overline Y) = \theta(X_u) = 0$ and 
$d\theta(\xi, \; \cdot\;) = 0$.  Since our CR structure is strictly
pseudoconvex,  equation~\eqref{eqn:Xu-equation} uniquely determines
$X_u$.   
Thus $H_0$ is the image of the first order differential
operator 
$\rho : \Gamma(M,F) \to \Gamma(M,T')$ defined by 
$\rho(u\xi) = X_u+u\xi$.

Define a second-order operator 
$D : \Gamma(M,F) \to \Gamma(M,E_1)$ 
as the composition $D = {\overline \partial_{T'}} \circ \rho$.  We then 
clearly get a complex 
\begin{equation}
\label{eqn:D-complex}
0 \longrightarrow \Gamma(M,F)
  \stackrel{D}{\longrightarrow} \Gamma(M,E_1)
  \stackrel{\overline{\partial}_1}{\longrightarrow} \Gamma(M,E_2)
  \stackrel{\overline{\partial}_2}{\longrightarrow} \cdots.
\end{equation}
It is this complex that we will use to get our subelliptic estimate
and therefore our decomposition theorems.  

Notice that $X_u$ includes a first derivative of $u$.  Using a local 
moving frame 
$\{e_1,\ldots,e_{n-1}\}$ for $^0T^{\prime}$
satisfying 
\begin{equation}
\label{eqn:local-frame}
L_{\theta}(e_{\alpha},\overline e_{\beta}) =  \delta_{\alpha\beta}, 
\end{equation}
we set $X_u = X^{\alpha} e_{\alpha}$ (note implicit sum).  Expanding 
$\theta([\overline{e}_{\beta},X^{\alpha} e_{\alpha} + u\xi]) = 0$, 
we get 
\begin{equation}
\theta\left(
	(\overline{e}_{\beta} X^{\alpha}) e_{\alpha} 
	+ X^{\alpha} [\overline{e}_{\beta},e_{\alpha}]
	+ (\overline{e}_{\beta} u) \xi + u [\overline{e}_{\beta},\xi]
      \right) = 0.
\end{equation}
This simplifies to 
$X^{\alpha} \left( - i \delta_{\beta\alpha} \right) + 
\overline{e}_{\beta} u = 0$,
so $X^{\alpha} = i \delta^{\alpha\beta}\overline e_{\beta} u$. 
Thus $\rho$ is indeed a first order operator, and our composition $D =
{\overline \partial_{T'}} \circ \rho$ is a second order operator.

Finally, we would like to relate our operator $D$ to that of Rumin
\cite{Rumin:1994}.  Define, for $p+q \ge n$, 
\begin{equation}
\label{eqn:Fpq-defn}
F^{p,q} = \{ u \in \theta \wedge \Lambda^{p-1,q} H^*
	: d\theta \wedge u = 0 \},
\end{equation}
and set $F^k = \oplus_{p+q=k} F^{p,q}$ for $k\ge n$.  Although the
definition \eqref{eqn:Fpq-defn} seems to depend on non-invariant
decomposition \eqref{eqn:bigrading}, we may actually express $F^k$
invariantly as 
\begin{equation*}
F^k = \{ u\in \C \otimes \Lambda^k M : 
	v \wedge u = 0 \text{ for all } 
	v \in \langle \theta, d\theta \rangle
      \},
\end{equation*}
where $\langle \theta, d\theta \rangle$ is the ideal generated by
$\theta$ and $d\theta$.  Since this ideal is CR-invariant, the
definition of $F^k$ is as well.  Below the middle dimension, we define
a slightly different space.  For $p+q=k \le n-1$, set 
\begin{equation*}
E^{p,q} = \Lambda^{p,q} H^* / \langle d\theta \rangle
\end{equation*}
and $E^{k} = \oplus_{p+q=k} E^{p,q}$, so 
\begin{equation*}
E^{k} = \C \otimes \Lambda^{k} M / \langle \theta, d\theta \rangle
\end{equation*}
is CR-invariant as well.   Rumin's $D$ operator is a map 
$D : E^{n-1} \to F^{n}$ given by $D[u] = d\tilde{u}$, where the
representative $\tilde{u}$ of  $[u] \in E^{n-1}$ is chosen so that
$d\tilde{u}$ will be in $F^{n}$.  There is then a complex 
\begin{equation}
\label{eqn:Rumin-k-complex}
\cdots \stackrel{d}{\longrightarrow}
E^{n-1} \stackrel{D}{\longrightarrow}
F^{n} \stackrel{d}{\longrightarrow}
F^{n+1} \stackrel{d}{\longrightarrow} \cdots,
\end{equation}
which decomposes into subcomplexes 
\begin{equation}
\label{eqn:Rumin-pq-complex}
\cdots \stackrel{d''}{\longrightarrow}
E^{p,n-p-1} \stackrel{D''}{\longrightarrow}
F^{p,n-p} \stackrel{d''}{\longrightarrow}
F^{p,n-p+1} \stackrel{d''}{\longrightarrow} \cdots.
\end{equation}
We hope to provide more details on these complexes in another paper.

The relation between our complex \eqref{eqn:D-complex} and Rumin's
complex \eqref{eqn:Rumin-pq-complex} occurs when $p=n-1$ in Rumin's
complex, in which case \eqref{eqn:Rumin-pq-complex} is 
\begin{equation}
\label{eqn:Rumin-complex-n-1}
0 \longrightarrow
E^{n-1,0} \stackrel{D''}{\longrightarrow}
F^{n-1,1} \stackrel{d''}{\longrightarrow}
F^{n-1,2} \stackrel{d''}{\longrightarrow} \cdots
\end{equation}
and we note that $E^{n-1,0} = \Lambda^{n-1,0} H^* = \Lambda^{n-1}
({}^0T^{\prime})^*$.  Let $K_M$ denote a nonvanishing closed $(n,0)$-form
(that is, an element of $\theta \wedge \Lambda^{n-1,0} H^*$), if one
exists.  For any positive $k$, we get a map $P_k : \Gamma(M,E_k) \to
F^{n-1,k}$ by interior multiplying the vector part of $u \in
\Gamma(M,E_k)$ into $K_M$, then wedging the remainder with the form
part of $u$.  Let $P_0 : \Gamma(M,F) \to E^{n-1,0}$ be given by $P_0
(u\xi) = u\cdot K_M$.  The claim is that each $P_k$ is an isomorphism
and the following diagram commutes:
\begin{align*}
\begin{CD}
0 @>>> \Gamma(M,F)    @>{D}>> 
	\Gamma(M,E_1) @>{\overline \partial_1}>> 
	\Gamma(M,E_2) @>{\overline \partial_2}>> 
	\cdots \\
  &&   @V{P_0}VV  @V{P_1}VV @V{P_2}VV  \\
0 @>>> E^{n-1,0} @>{D''}>> F^{n-1,1} @>{d''}>> F^{n-1,2} @>{d''}>> \cdots \\
\end{CD}
\end{align*}
Since $K_M$ always exists locally, the two complexes are locally
isomorphic.  If the canonical line bundle is trivial, then this complex
version \eqref{eqn:Rumin-pq-complex} of the Rumin complex is isomorphic 
to our new complex \eqref{eqn:D-complex}.

\section{A subelliptic estimate and decomposition theorem}
In this section, we state two of our main results.  First, we produce
a subelliptic estimate at $\Gamma(M,E_1)$ for our
complex~\eqref{eqn:D-complex} in the $5$-dimensional case.  
Using this, we get a
Hodge-Kodaira decomposition theorem for elements of $\Gamma(M,E_1)$.  

We begin with some preliminaries.  
Our choice of pseudohermitian structure $\theta$ 
determines the \emph{pseudohermitian connection} $\nabla$ (see
\cite{Webster:1978,Tanaka:book}): this is the unique connection that
is compatible with $H$ and its complex structure, for which $\theta$ and
$d\theta$ are parallel, and satisfying an additional
torsion condition.  For any tensor field $u$ on $M$, the total
covariant derivative $\nabla u$ can be decomposed as
\begin{equation*}
\nabla u = \nabla' u + \nabla'' u + \nabla_T u \otimes \theta,
\end{equation*}
where $\nabla' u$ involves derivatives only with respect to vector
fields in ${}^0T^{\prime}$, and $\nabla'' u$ only with respect to
vector fields in ${}^0T^{\prime\prime}$.  Writing $\nabla_H u =
\nabla'u + \nabla''u$, the \emph{Folland-Stein norms} $\norm{\cdot}_k$
are defined by 
\begin{equation*}
\norm{u}_k^2 = \sum_{j=0}^k \norm{\nabla_H^j u}^2,
\end{equation*}
where $\norm{\cdot}$ denotes the $L^2$ norm defined with respect
to the Webster metric.  (Note that in
\cite{Akahori:1981}, the $\norm{\cdot}_1$ and $\norm{\cdot}_2$ norms
were called $\norm{\cdot}^{\prime}$ and $\norm{\cdot}^{\prime\prime}$,
respectively.)  We will write $(\ ,\ )$ for the hermitian inner
product that corresponds to the norm $\norm{\cdot}$, and for 
any bundle $E$ we will let $\Gamma_2(M,E)$
denote the completion of $\Gamma(M,E)$ with respect to the 
$L^2$ norm.

Define a second-order operator 
$L = 1 + {\nabla^{\prime}}^{\ast} \nabla^{\prime}
	+ {\nabla^{\prime\prime}}^{\ast} \nabla^{\prime\prime}$.  We
then define our Laplacian $\square : \Gamma(M,E_1) \to \Gamma(M,E_1)$ by 
$\square u = D D^* u + \overline \partial_1^* L \overline \partial_1 u$,
where the adjoints are defined with respect to the
complex~\eqref{eqn:D-complex}.  
We use this operator and the norms defined above to express our
subelliptic estimate in the following theorem.
\begin{Theorem}[Main Estimate]
\label{thm:main-estimate}
Let $(M,{}^0T^{\prime\prime})$ be a compact, strictly pseudoconvex CR
manifold of dimension $5$.  Then there exists a constant $c > 0$ such
that  
\begin{equation}
\label{eqn:main-estimate}
\left( \phi, \square \phi \right) 
  = \norm{ D^{\ast}\phi }^2 
	+ \norm{ \overline \partial_1 \phi }_1^2
   \geq c \norm{ \phi }_2^2 - \norm{ \phi }_1^2
\end{equation}
for all $\phi \in \Gamma(M,E_1)$.
\end{Theorem}
The details of the proof of this estimate will be confined to
the next section.

We define new norms that are Sobolev extensions of the
Folland-Stein norms $\norm{\cdot}_k$ as follows.  We set 
\begin{equation*}
\norm{ u }_{k,m}^2 
	=  \sum_{l = 0}^{m} 
	   \sum_{j=0}^k 
		\norm{ \nabla^{l} \nabla_H^j u }^2.
\end{equation*}
The first parameter, $k$, specifies the number of derivatives in the
$H$ directions, whereas the second parameter, $m$, is the number of
unconstrained derivatives.  (We remark that in \cite{Akahori:1981}
these norms were written slightly differently: for example,
$\norm{\cdot}_{2,m}$ was $\norm{\cdot}^{\prime\prime}_{(m)}$.)  Then
our main estimate, Theorem~\ref{thm:main-estimate}, together with
standard integration-by-parts techniques, gives us the following
Sobolev estimate.
\begin{corollary}
\label{thm:main-Sobolev-estimate}
Let $(M,{}^0T^{\prime\prime})$ be a compact, strictly pseudoconvex CR
manifold of dimension $5$.  For each positive integer $m$, there
exists a constant $c_m > 0$ such that  
\begin{equation}
\label{eqn:main-Sobolev-estimate}
\norm{ D^{\ast}\phi }^2_{0,m} 
	+ \norm{ \overline \partial_1 \phi }_{1,m}^2
   \geq c_m \norm{ \phi }_{2,m}^2 - \norm{ \phi }_{1,m}^2
\end{equation}
for all $\phi \in \Gamma(M,E_1)$.
\end{corollary}

Let us write $\Harmonics$ for the harmonic elements of
$\Gamma(M,E_1)$, with respect to the Laplacian $\square$.  
In order to find a useful expression for $\Harmonics$, we use the
following lemma to express the adjoint of $D$ in simpler terms.
\begin{lemma}
Let $\widetilde{H}_0$ be the completion of $H_0$ under the $L^2$
norm, and $\pi_{\widetilde{H}_0} : \Gamma_2(M,T') \to \widetilde{H}_0$
is orthogonal projection.  Then we have the following relations:
\begin{enumerate}
\renewcommand{\labelenumi}{(\alph{enumi})}
\item
$\overline \partial_0^* = \pi_{\widetilde{H}_0} \circ \overline
\partial_{T'}^*$, where $\overline \partial_{T'}^*$ is the formal adjoint of
$\overline \partial_{T'}$

\item
$\ker D^*  = \ker \overline \partial_0^*$

\end{enumerate}
\end{lemma}
\begin{proof}
The first conclusion follows from the relation between the standard
deformation complex~\eqref{eqn:standard-complex} and the
complex~\eqref{eqn:H0-complex} involving $H_0$.   Since 
$H_0 \subset \Gamma(M,T')$ and 
$\Gamma(M,E_1) \subset \Gamma(M,T'\otimes ({}^0T'')^*)$, we may write
$\overline\partial_0 = \overline \partial_{T'} \circ \pi_{\widetilde{H}_0}$, 
from which it follows that 
$\overline \partial_0^* = \pi_{\widetilde{H}_0} \circ \overline
\partial_{T'}^*$ on $\Gamma(M,E_1)$.  
That $\ker D^* = \ker \overline \partial_0^*$ is due to two simple
facts: first, that $D^* = \rho^* \circ \overline\partial_0^*$ and,
second, that $\rho : \Gamma(M,F) \to H_0$ is an isomorphism. 
\end{proof}

This lemma then implies that we may write $\Harmonics$ as 
\begin{equation*}
\Harmonics
	= \ker \square 
	= \{ \phi \in \Gamma(M,E_1) : 
		\overline \partial_0^* \phi = 0 
		\quad \text{and} \quad 
		\overline \partial_1 \phi = 0 
	  \}.
\end{equation*}
The subelliptic estimate in Theorem~\ref{thm:main-estimate} gives us
the following Hodge-Kodaira decomposition theorem.
\begin{Theorem} 
\label{thm:Hodge-Kodaira}
Let $(M,{}^0T^{\prime\prime})$ be a compact, strictly pseudoconvex CR
manifold of dimension $5$.  Then
\begin{equation*}
\Harmonics \cong \frac{ \ker \overline \partial_1 }{ \im D }.
\end{equation*}
Moreover, there exists a Neumann operator 
$N : \Gamma_2(M,E_1) \to \Gamma_2(M,E_1)$ and a harmonic projector
$H : \Gamma_2(M,E_1) \to \mathcal{H}$
satisfying 
$NH = HN = 0$,
$[N,DD^*] = 0 = [N,\overline \partial_1^* L \overline \partial_1]$,
and 
$u= Hu + \square Nu = Hu + N\square u $
for all $u\in \Gamma_2(M,E_1)$.
\end{Theorem}

We will construct the Neumann operator $N$ and the harmonic projector $H$
by considering the differential equation 
\begin{equation}
\label{eqn:our-Poisson}
\square u = f.
\end{equation}
Let us write $\Harmonics^{\perp}$ for elements of $\Gamma_2(M,E_1)$
that are orthogonal to $\Harmonics$ with respect to the $L^2$ norm.
We begin with a fairly standard lemma.
\begin{lemma}
\label{lemma:Prop-8.1}
There is a constant $c>0$ for which
\begin{equation*}
\norm{ D^* u }^2+\norm{ \overline \partial_1 u }^2_1 \ge c \norm{u}_1^2
\end{equation*}
for all $u \in \Harmonics^{\perp} \subset \Gamma_2(M,E_1)$.
\end{lemma}
\begin{proof}[Proof of Lemma~\ref{lemma:Prop-8.1}]
We assume the conclusion is false.  That is, for each integer $k>0$,
we assume that there is an element $u_k \in \Harmonics^{\perp}$ satisfying  
$\norm{ D^* u_k }^2+\norm{ \overline \partial_1 u_k }^2_1 \le \frac{1}{k} \norm{ u_k }_1^2$.
Rescaling these $u_k$ if necessary, we may assume that 
$\norm{ u_k }_1 = 1$ and therefore 
$\norm{ D^* u_k }^2+\norm{ \overline \partial_1 u_k }^2_1 \le \frac{1}{k}$.
By our estimate~\eqref{eqn:main-estimate} (extended by continuity
to $\Gamma_2(M,E_1)$), we have
\begin{align*}
c\norm{ u_k }^2_2
	& \le    
			\norm{ D^* u_k }^2
			+ \norm{ \overline \partial_1 u_k }^2_1
			+ \norm{ u_k }^2_1
		  \\
	& \le   \left( \frac{1}{k} + 1 \right) \le 2.
\end{align*}
The sequence $\{u_k\}$ is thus bounded with
respect to $\norm{\cdot}_{2}$, the Folland-Stein $2$-norm.  Any such
set is precompact with respect to $\norm{\cdot}_{1}$; this means
there is a subsequence $\{ u_{k_j} \}$ that converges weakly 
in $\Gamma_2(M,E_1)$ and strongly in the 
Folland-Stein $1$-norm.  Let $u$ be its limit.  On the one hand, $u
\in \Harmonics^{\perp}$ as each element $u_{k_j}$ is.  On the other
hand, the closedness of the differential operator $\square$ implies that
$u \in \Dom \square$ and $\square u = 0$.  Thus $u \in \Harmonics$, so $u =
0$.  But $\norm{u}_1 = 1$, so this is a contradiction.  
\end{proof}

\begin{proof}[Proof of Theorem \ref{thm:Hodge-Kodaira}]
By
Lemma~\ref{lemma:Prop-8.1} and Theorem~\ref{thm:main-estimate}, 
the quadratic form
\begin{equation*}
Q(u,u) =    
			\norm{ D^* u }^2
			+ \norm{ \overline \partial_1 u }^2_1
\end{equation*}
defines a norm that
is equivalent to $\norm{\cdot}_{2}$.  We endow
$\Harmonics^{\perp}$ with this norm, and let $Q(u,v)$ denote the
associated symmetric bilinear form.  Note that if $u$ and $v$ are
smooth, then $Q(u,v) = (\square u,v)$.

By
Lemma~\ref{lemma:Prop-8.1}, the linear functional $v \mapsto (f,v)$ is
bounded on
$\Harmonics^{\perp}$ for any $f \in \Gamma_2(M,E_1)$.  The Riesz
representation theorem then implies that there is a unique $u \in
\Harmonics^{\perp}$ such that $Q(  u, v ) = \left( f,
v\right)$ for all $v \in \Harmonics^{\perp}$.  Thus we have solved
\eqref{eqn:our-Poisson} for $f \in \Harmonics^{\perp}$.

The Neumann operator is given by $Nf = u$, the solution 
$u \in \Harmonics^{\perp}$ to $\square u =
f$ in the above sense.  This makes sense for $f \in
\Harmonics^{\perp}$, so under the orthogonal decomposition 
$\Gamma_2(M,E_1) = \Harmonics \oplus \Harmonics^{\perp}$
we can extend $N$ to all of $\Gamma_2(M,E_1)$ by declaring that it is
identically zero on $\Harmonics$.  We define the harmonic projector
$H$ as orthogonal projection onto $\Harmonics$ under this
decomposition.  The operators $H$ and $N$ project onto orthogonal
spaces, so $HN = 0 = NH$.  On the other hand, the decompositions
$u = Hu + \square N u = Hu + N \square u$ follow immediately from the
construction of $N$ and $H$.  

To see that 
$[\overline \partial_1^* L \overline \partial_1, N] = 0
 = [DD^*, N]$
takes a bit more work.  From $[\square,N]=0$ it follows directly
that 
$[\overline \partial_1^* L \overline \partial_1, N] + [DD^*, N] = 0$,
so we need only show that, say, $[DD^*, N] = 0$.  
This follows easily by considering separately $u\in \Harmonics$
(on which $DD^*$ and $N$ are separately zero) and 
$u=\square v\in \Harmonics^{\perp}$,
in which case $[DD^*, N]\square v = 0$ is a straightforward
computation based on the formulas
$N\square v = v - Hv$, $[DD^*, \square] = 0$, and 
$HDD^* = DD^*H = 0$.

%
%

Finally, the isomorphism $\Harmonics \cong \ker \overline\partial_1/
\im D$ follows as usual from the existence of the Neumann
operator, since the harmonic projector $H$
restricts to a map $H\colon \ker\overline\partial_1 \to \Harmonics$
whose kernel is exactly $\im D$ by the arguments above.
\end{proof}

\section{Proof of the subelliptic estimate}
In this section, we prove Theorem~\ref{thm:main-estimate}, our
subelliptic estimate.  
Since our manifold $M$ is assumed to be compact, it will suffice to show
that \eqref{eqn:main-estimate} holds for $\phi$ supported in
a neighborhood of each point: assuming this, we can choose
a locally finite collection $\{\alpha_i\}$ of smooth nonnegative functions
satisfying $\sum_i\alpha_i^2=1$, apply 
\eqref{eqn:main-estimate} to $\alpha_i\phi$, and sum over $i$,
yielding \eqref{eqn:main-estimate} plus some lower-order terms
that can be absorbed into the right-hand side.

Let $\{e_1,e_2\}$ be a local moving frame for ${}^0T'$ satisfying
\eqref{eqn:local-frame}, from which it follows that 
\begin{equation}
\label{eqn:frame-comm}
\pi_F [e_{\alpha}, \overline{e}_{\beta}] = - i \delta_{\alpha\beta} \ \xi, 
\end{equation}
and let $\{\theta^1, \theta^2\}$ be the 
dual sections of $({}^0T^{\prime})^*$, thought of as
one-forms according to the decomposition \eqref{eqn:vector-bundles-2}.  
We may then write $\phi\in 
\Gamma(M,{}^0T^{\prime} \otimes \Lambda^j(^0T^{\prime\prime})^{\ast})$
in coordinates as
\begin{equation}
\phi
	 = \phi^{\alpha}_{\beta_1,\ldots,\beta_j}
	   e_{\alpha}
	   \otimes \overline \theta^{\beta_1} \wedge \cdots \wedge
	   \overline \theta^{\beta_j}.
\end{equation}
(Notice the implicit sums over $\alpha$ and $\beta_1$ through
$\beta_j$.)  Throughout this section,
we will assume $\phi$ is supported in the neighborhood on which our
moving frame is defined, so that 
\begin{equation}
\label{eqn:local-norm-phi}
\norm{\phi}^2 
	= \sum_{\alpha,\beta_1,\ldots,\beta_j} 
		\norm{\phi^{\alpha}_{\beta_1,\ldots,\beta_j} }^2 .
\end{equation}

We will often find it useful to look only at the top order
derivatives.  In light of the commutation
relation~\eqref{eqn:frame-comm}, this unfortunately is not possible.
Instead, we will look at only the top \emph{weight} derivatives, where
we allocate a weight of $1$ to vector fields in $H$ and a weight of $2$
to $\xi$.  We will then write $\sim$ for
equal modulo lower weight terms.  This generalizes to
$\gtrsim$ and $\lesssim$, meaning greater than or less than, modulo
negligible terms.  Our main estimate~\eqref{eqn:main-estimate} can
thus be written  
\begin{equation*}
\left( u, \square u \right) 
  = \norm{D^{\ast}\phi}^2 
	+ \norm{\overline \partial_1 \phi }_1^2
   \gtrsim c \norm{\phi}_2^2
\end{equation*}
for all $\phi \in \Gamma(M,E_1)$.
To prove this estimate, we will need a local expression for
$\norm{\phi}_2^2$ rather than $\norm{\phi}^2$.  Modulo lower weight
terms, this expression is
\begin{equation*}
\norm{ \phi }^2_2 
  \sim
  \sum_{k,l,\alpha,\beta_1,\ldots,\beta_j} 
	\norm{ 
	  e_k e_{l} \phi^{\alpha}_{\beta_1,\ldots,\beta_j}
 	}^2,
\end{equation*}
where $j,k \in \{1,\ldots,n-1,\overline 1, \ldots,\overline{n-1} \}$.

We begin the actual proof of Theorem~\ref{thm:main-estimate} by
describing $\phi$, $\overline \partial_1 \phi$, and $D^* \phi$ in terms
of our local moving frame.   (Compare 
\cite[Lemmas 3.2 and 3.3, page 317]{Akahori:1981}.)
\begin{lemma}
\label{lemma:phi-lemma}
Suppose $\phi \in \Gamma(M,E_1)$.  Then $\phi^{1}_{2} = \phi^{2}_{1}$,
$(\overline \partial_1 \phi)^{\alpha}_{1,2}
 \sim \overline e_1 \phi^{\alpha}_2 - \overline e_2 \phi^{\alpha}_1$,
and
\begin{equation}
\label{eqn:D*phi-frame} 
D^* \phi \sim - i \left( 
	  e_1 e_1\phi^1_1 
	+ e_1 e_2\phi^1_2 
	+ e_2 e_1\phi^2_1 
	+ e_2 e_2\phi^2_2
		  \right)  \xi.
\end{equation}
\end{lemma}
\begin{proof}
In our local frame, we may write $\phi = \phi^{\alpha}_{\beta}
e_{\alpha} \otimes \overline \theta^{\beta}$.  (Since $\Gamma(M,E_1)
\subset \Gamma(M,{}^0T^{\prime} \otimes ({}^0T^{\prime\prime})^*)$,
there are no $\xi \otimes \overline \theta^{\beta}$ terms.)  In this
case 
$\overline \partial^{(1)}
\phi (\overline e_2, \overline e_2)$ 
is (see equation~\eqref{eqn:delbarphi-defn})
\begin{align}
\overline \partial^{(1)} \phi (\overline e_1, \overline e_2)
 & =   \pi' [\overline e_1, \phi (\overline e_2) ]
	- \pi' [\overline e_2, \phi (\overline e_1) ]
	- \phi ([\overline e_1, \overline e_2]) 
	\nonumber \\
 & =   (\overline e_1 \phi^{\alpha}_2) e_{\alpha}
        + \phi^{\alpha}_2 \pi' [\overline e_1, e_{\alpha}] 
	- (\overline e_2 \phi^{\alpha}_1 ) e_{\alpha} 
	\label{eqn:delbarphi-lemma}
	\\
 &      - \phi^{\alpha}_1 \pi' [\overline e_2, e_{\alpha}]
	- \phi ([\overline e_1, \overline e_2]) 
	\nonumber \\
 & \sim  
	(\overline e_1 \phi^{\alpha}_2) e_{\alpha}
       	- (\overline e_2 \phi^{\alpha}_1 ) e_{\alpha},
	\nonumber
\end{align}
where we have discarded all the terms without a derivative of a
component of $\phi$.  This proves the second claim; the first claim
follows from applying the one-form $\theta$ to both sides of
\eqref{eqn:delbarphi-lemma}: 
\begin{equation*}
0 = \phi^{\alpha}_2 \theta([\overline e_1, e_{\alpha}])
    - \phi^{\alpha}_1 \theta([\overline e_2, e_{\alpha}])
  = i \phi^{\alpha}_2 \delta_{1\alpha}
    - i \phi^{\alpha}_1 \delta_{2\alpha}
  = i (\phi^{1}_2 - \phi^2_1),
\end{equation*}
where we have simplified using
$\theta([e_{\alpha}, \overline e_{\beta}]) = -i \delta_{\alpha \beta}$.

Finally, we prove equation~\eqref{eqn:D*phi-frame}.  To compute this
adjoint, we take the inner product of $D^*\phi$ with an element $u\xi$
of $\Gamma(M,F)$, and integrate by parts:
\begin{equation*}
(u\xi, D^*\phi) = (D(u\xi), \phi) \sim 
	\left( \overline \partial_0 
		(u\xi + i (\overline e_1 u) e_1 
			+ i (\overline e_2 u) e_2),
		\phi \right).
\end{equation*}
If we write $\psi = \psi^{\alpha}_{\overline \beta} \ e_{\alpha}
\otimes \overline \theta^{\beta}$ for 
$D (u\xi) = \overline \partial_0 (u\xi + i (\overline e_1 u) e_1 + i
(\overline e_2 u) e_2)$
(again, there is no $\xi \otimes \overline \theta^{\beta}$ term as
$D(u\xi) \in \Gamma(M,E_1)$),
then we can compute
$\psi^{\alpha}_{\overline \beta} 
	= \theta^{\alpha} ( \psi( \overline e_{\beta} ) )$.
The inside term is not difficult to compute, and we get 
$\psi( \overline e_{\beta} ) \sim 
	\pi' \left[ \overline e_{\beta},
		u\xi + i (\overline e_1 u) e_1 
		     + i (\overline e_2 u) e_2
	     \right]$,
so $\psi^{\alpha}_{\overline \beta} \sim 
	i \overline e_{\beta} \overline e_{\alpha} u$. 
Undoing the integration by parts above gives 
equation~\eqref{eqn:D*phi-frame}. 
%
%
\end{proof}

The primary tool in our proof of Theorem~\ref{thm:main-estimate} is
the following lemma.  This follows at least in part from the local
expressions computed in Lemma~\ref{lemma:phi-lemma}
\begin{lemma}[Key Estimate]
\label{lemma:key-estimate}
For all $\phi \in \Gamma(M,E_1)$, 
\begin{eqnarray}
\lefteqn{
    \norm{ D^* \phi }^2 + 2 \norm{ \overline \partial_1 \phi }^2_1
    } \nonumber \\
  & \gtrsim & 
	\norm{ e_1 e_1\phi^1_1 }^2 + \norm{ e_2 e_2\phi^2_2 }^2 
	+ 4 \norm{ e_1 \overline e_2 \phi^1_2 }^2 
	\label{eqn:key-estimate}
	\\
  & &   
	+ 4 \norm{ e_2 \overline e_1 \phi^1_2 }^2
	+ \norm{ \overline e_1 \overline e_1 \phi^2_2 }^2 
	+ \norm{ \overline e_2 \overline e_2 \phi^1_1 }^2  
	\nonumber
\end{eqnarray}
\end{lemma}
\begin{proof}
We begin by computing $\norm{ D^* \phi }^2$.  From
\eqref{eqn:D*phi-frame}, we have 
\begin{equation*}
\norm{ D^* \phi }^2  \sim 
	\norm{ e_1 e_1\phi^1_1 
	+ e_1 e_2\phi^1_2 
	+ e_2 e_1\phi^2_1 
	+ e_2 e_2\phi^2_2 
	}^2.
\end{equation*}
We expand this to get 
\begin{align}
\norm{ D^* \phi }^2 
    & \sim \norm{ e_1 e_1\phi^1_1 }^2 + \norm{ e_1 e_2\phi^1_2 }^2 
	   + \norm{ e_2 e_1\phi^2_1 }^2 + \norm{ e_2 e_2\phi^2_2 }^2 
	   \nonumber \\
  &     \quad   + 2 \Real (e_1 e_1\phi^1_1, e_1 e_2\phi^1_2)
	   + 2 \Real (e_1 e_1\phi^1_1, e_2 e_1\phi^2_1) 
	   \nonumber \\ 
  &     \quad   + 2 \Real (e_1 e_1\phi^1_1, e_2 e_2\phi^2_2) 
	   + 2 \Real (e_1 e_2\phi^1_2, e_2 e_1\phi^2_1) 
	   \label{eqn:D*phi-1}
	   \\
  &   \quad     + 2 \Real (e_1 e_2\phi^1_2, e_2 e_2\phi^2_2)
  	   + 2 \Real (e_2 e_1\phi^2_1, e_2 e_2\phi^2_2).
	   \nonumber
\end{align}
Since $\phi^{1}_2 = \phi^2_1$ by Lemma~\ref{lemma:phi-lemma} and
$[e_{\alpha}, e_{\beta}] \sim 0$ for all $\alpha$ and $\beta$, one of
the cross terms simplifies:  
$2 \Real (e_1 e_2\phi^1_2, e_2 e_1\phi^2_1) 
	= 2 \Real (e_1 e_2\phi^1_2, e_1 e_2\phi^1_2) 
	= 2 \norm{ e_1 e_2\phi^1_2 }^2$.
Four of the other cross terms combine and \eqref{eqn:D*phi-1}
simplifies to  
\begin{align}
\norm{ D^* \phi }^2 
  &\sim  
	\norm{ e_1 e_1\phi^1_1 }^2 + \norm{ e_1 e_2\phi^1_2 }^2 
	+ \norm{ e_2 e_1\phi^1_2 }^2 + \norm{ e_2 e_2\phi^2_2 }^2 
	\nonumber \\
  &  \quad
	+ 4 \Real (e_1 e_1\phi^1_1, e_1 e_2\phi^1_2)
	+ 2 \Real (e_1 e_1\phi^1_1, e_2 e_2\phi^2_2) 
	\label{eqn:D*phi-2}
	\\ 
  &  \quad
	+ 2 \norm{ e_1 e_2\phi^1_2 }^2
	+ 4 \Real (e_1 e_2\phi^1_2, e_2 e_2\phi^2_2).
	\nonumber 
\end{align}

We will deal with the remaining cross terms by adding 
$2 \norm{ \overline \partial_1 \phi }^2_1$.  By Lemma \ref{lemma:phi-lemma},  
\begin{align*}
2 \norm{ \overline \partial_1 \phi }^2_1
  &\sim  
	2 \norm{ \overline e_1 \phi^1_2 - \overline e_2 \phi^1_1 }_1^2 
	+ 2 \norm{ \overline e_1 \phi^2_2 - \overline e_2 \phi^1_2 }_1^2
  \\
  &\sim  
	2 \norm{ e_1 (\overline e_1 \phi^1_2 - \overline e_2 \phi^1_1) }^2 
	+ 2 \norm{ e_2(\overline e_1 \phi^1_2 - \overline e_2 \phi^1_1) }^2  \\
   &  \quad  
	+ 2 \norm{ \overline e_1 
		(\overline e_1 \phi^1_2 - \overline e_2 \phi^1_1) }^2 
	+ 2 \norm{ \overline e_2
		(\overline e_1 \phi^1_2 - \overline e_2 \phi^1_1) }^2  \\
   &   \quad
	+ 2 \norm{ e_1 (\overline e_1 \phi^2_2 - \overline e_2 \phi^1_2) }^2
	+ 2 \norm{ e_2 (\overline e_1 \phi^2_2 - \overline e_2 \phi^1_2) }^2 \\
   &   \quad
	+ 2 \norm{ \overline e_1 
		(\overline e_1 \phi^2_2 - \overline e_2 \phi^1_2) }^2
	+ 2 \norm{ \overline e_2
		(\overline e_1 \phi^2_2 - \overline e_2 \phi^1_2) }^2
  \\
  &\gtrsim   
	2 \norm{ e_1 (\overline e_1 \phi^1_2 - \overline e_2 \phi^1_1) }^2 
	+ 2 \norm{ e_2 (\overline e_1 \phi^2_2 - \overline e_2 \phi^1_2) }^2 \\
  &\quad  
	+ 2 \norm{ \overline e_2
		(\overline e_1 \phi^1_2 - \overline e_2 \phi^1_1) }^2  
	+ 2 \norm{ \overline e_1 
		(\overline e_1 \phi^2_2 - \overline e_2 \phi^1_2) }^2
\end{align*}
Since $[e_{\alpha},\overline{e_{\beta}}] \sim -i \delta_{\alpha\beta}
\xi$, we have
\begin{align*}
e_1 (\overline e_1 \phi^1_2 - \overline e_2 \phi^1_1) 
   & \sim  
	-i \xi \phi^1_2 + \overline e_1 e_1 \phi^1_2 
	- e_1 \overline e_2 \phi^1_1 \\
e_2 (\overline e_1 \phi^2_2 - \overline e_2 \phi^1_2) 
   & \sim 
	e_2 \overline e_1 \phi^2_2 - \overline e_2 e_2 \phi^1_2 
	+ i \xi \phi^1_2.
\end{align*}
Moreover, 
\begin{align*}
2 \norm{ \overline e_2
	(\overline e_1 \phi^1_2 - \overline e_2 \phi^1_1) }^2  
&+ 2 \norm{ \overline e_1 
	(\overline e_1 \phi^2_2 - \overline e_2 \phi^1_2) }^2
 \\
& \ge 
\norm{ \overline e_2
	(\overline e_1 \phi^1_2 - \overline e_2 \phi^1_1) 
   + \overline e_1 
	(\overline e_1 \phi^2_2 - \overline e_2 \phi^1_2) }^2 \\
& \sim   
	\norm{ \overline e_1 \overline e_1 \phi^2_2
	   - \overline e_2 \overline e_2 \phi^1_1  }^2
\end{align*}
as $[\overline e_{\alpha}, \overline e_{\beta}] \sim 0$.
Hence 
\begin{align}
2 \norm{ \overline \partial_1 \phi }^2_1
  &\gtrsim 
	2 \norm{ 
		-i \xi \phi^1_2 + \overline e_1 e_1 \phi^1_2 
		- e_1\overline e_2 \phi^1_1 
	  }^2 \nonumber \\
  &  \quad
	+ 2 \norm{ 
		e_2 \overline e_1 \phi^2_2 - \overline e_2 e_2 \phi^1_2 
		+ i \xi \phi^1_2 
	  }^2 \nonumber \\
  &  \quad
	+ \norm{ 
		\overline e_1 \overline e_1 \phi^2_2
		- \overline e_2 \overline e_2 \phi^1_1  
	  }^2 \nonumber \\
  &\sim 
	2 \norm{ \xi \phi^1_2 }^2 
	+ 2 \norm{ \overline e_1 e_1 \phi^1_2 }^2 
	+ 2 \norm{ e_1 \overline e_2 \phi^1_1 }^2 \nonumber \\
  & \quad
	-4 \Real ( i \xi \phi^1_2, \overline e_1 e_1 \phi^1_2 )
	+4 \Real ( i \xi \phi^1_2, e_1 \overline e_2 \phi^1_1 )
	\nonumber \\
  &  \quad
	-4 \Real ( \overline e_1 e_1 \phi^1_2, e_1 \overline e_2 \phi^1_1 ) 
	+ 2 \norm{ e_2 \overline e_1 \phi^2_2 }^2 
	\nonumber \\
  &  \quad
	+ 2 \norm{ \overline e_2 e_2 \phi^1_2 }^2
	+ 2 \norm{ \xi \phi^1_2 }^2  
	-4 \Real (e_2 \overline e_1 \phi^2_2, \overline e_2 e_2 \phi^1_2 )
	\label{eqn:delbarphi-norm}
	\\
  &  \quad
	+4 \Real (e_2 \overline e_1 \phi^2_2, i \xi \phi^1_2 )
	-4 \Real (\overline e_2 e_2 \phi^1_2, i \xi \phi^1_2 ) 
	\nonumber \\
  &  \quad
	+ \norm{ \overline e_1 \overline e_1 \phi^2_2 }^2 
	+ \norm{ \overline e_2 \overline e_2 \phi^1_1 }^2 
	-2 \Real (\overline e_1 \overline e_1 \phi^2_2,
		\overline e_2 \overline e_2 \phi^1_1 ).
	\nonumber
\end{align}

To cancel the cross terms in \eqref{eqn:D*phi-2}, 
we will make use of the fact that $e_2$ commutes with
$\overline e_1$ and $e_1$ modulo lower-weight terms,
and therefore integrating by parts yields
\begin{align*}
-4 \Real (e_2 \overline e_1 \phi^2_2, \overline e_2 e_2 \phi^1_2 )
 & \sim 
  - 4 \Real (\overline e_1 e_2 \phi^2_2, \overline e_2 e_2 \phi^1_2 )\\
 & \sim 
   4 \Real (e_2 \overline e_1 e_2 \phi^2_2, e_2 \phi^1_2 )\\
 & \sim 
   4 \Real (\overline e_1 e_2 e_2 \phi^2_2, e_2 \phi^1_2 )\\
 & \sim 
   -4 \Real (e_2 e_2 \phi^2_2, e_1 e_2 \phi^1_2 )\\
 & \sim 
  -4 \Real (e_1 e_2 \phi^1_2 , e_2 e_2 \phi^2_2).
\end{align*}
A similar argument shows that three of the cross terms on the
right-hand side of \eqref{eqn:delbarphi-norm} cancel 
all the cross terms of
\eqref{eqn:D*phi-2}: 
\begin{align*}
\norm{ D^* \phi }^2 + 2 \norm{ \overline \partial_1 \phi }^2_1
  &\gtrsim  
	\norm{ e_1 e_1\phi^1_1 }^2 + \norm{ e_1 e_2\phi^1_2 }^2 
	+ \norm{ e_2 e_1\phi^1_2 }^2 + \norm{ e_2 e_2\phi^2_2 }^2 
	\\
  &\quad 
	+ 2 \norm{ e_1 e_2\phi^1_2 }^2
	+ 2 \norm{ \xi \phi^1_2 }^2 
	+ 2 \norm{ \overline e_1 e_1 \phi^1_2 }^2 
	+ 2 \norm{ e_1 \overline e_2 \phi^1_1 }^2 
	\\
  &\quad  
	-4 \Real ( i \xi \phi^1_2, \overline e_1 e_1 \phi^1_2 )
	+4 \Real ( i \xi \phi^1_2, e_1 \overline e_2 \phi^1_1 )
	\\
  &\quad  
	+ 2 \norm{ e_2 \overline e_1 \phi^2_2 }^2 
	+ 2 \norm{ \overline e_2 e_2 \phi^1_2 }^2
	+ 2 \norm{ \xi \phi^1_2 }^2  
	\\
  &\quad  
	+4 \Real (e_2 \overline e_1 \phi^2_2, i \xi \phi^1_2 )
	-4 \Real (\overline e_2 e_2 \phi^1_2, i \xi \phi^1_2 ) 
	\\
  &  \quad
	+ \norm{ \overline e_1 \overline e_1 \phi^2_2 }^2 
	+ \norm{ \overline e_2 \overline e_2 \phi^1_1 }^2. 
\end{align*}
We now have more cross terms, this time involving $\xi$.  

We will deal with some of these cross terms using integration by parts.
The adjoint of $e_{\alpha}$ is $-\overline e_{\alpha}$, and so (using
$[e_2, \overline e_2] \sim - i \xi$ and other commutation relations),
we have 
\begin{align*}
-4 \Real ( i \xi \phi^1_2, \overline e_1 e_1 \phi^1_2 )
  & \sim 
	+ 4 \Real ( e_2 \overline e_2 \phi^1_2, 
		  \overline e_1 e_1 \phi^1_2 )
	- 4 \Real ( \overline e_2 e_2 \phi^1_2, 
		  \overline e_1 e_1 \phi^1_2 ) \\
  & \sim  
	- 4 \Real ( \overline e_2 \phi^1_2, 
		  \overline e_2 \overline e_1 e_1 \phi^1_2 )
	- 4 \Real ( \overline e_2 e_2 \phi^1_2, 
		  \overline e_1 e_1 \phi^1_2 ) \\
  & \sim  
	- 4 \Real ( \overline e_2 \phi^1_2, 
		  \overline e_1 e_1 \overline e_2 \phi^1_2 )
	- 4 \Real ( \overline e_2 e_2 \phi^1_2, 
		  \overline e_1 e_1 \phi^1_2 ) \\
  & \sim  
	+ 4 \norm{ e_1 \overline e_2 \phi^1_2 }^2 
	- 4 \Real ( \overline e_2 e_2 \phi^1_2, 
		  \overline e_1 e_1 \phi^1_2 ). 
\end{align*}
Similarly,
\begin{align*}
-4 \Real ( i \xi \phi^1_2, \overline e_2 e_2 \phi^1_2 )
  & \sim 
	+ 4 \Real ( e_1 \overline e_1 \phi^1_2, 
		  \overline e_2 e_2 \phi^1_2 )
	- 4 \Real ( \overline e_1 e_1 \phi^1_2, 
		  \overline e_2 e_2 \phi^1_2 ) \\
  & \sim  
	+ 4 \norm{ e_2 \overline e_1 \phi^1_2 }^2 
	- 4 \norm{ e_2 e_1 \phi^1_2 }^2.
\end{align*}
Thus 
\begin{align}
	\norm{ D^* \phi }^2& + 2 \norm{ \overline \partial_1 \phi }^2_1
 \nonumber \\
  &\gtrsim  
	\norm{ e_1 e_1\phi^1_1 }^2 + \norm{ e_1 e_2\phi^1_2 }^2 
	+ \norm{ e_2 e_1\phi^1_2 }^2 + \norm{ e_2 e_2\phi^2_2 }^2 
	\nonumber \\
  &\quad
	+ 2 \norm{ e_1 e_2\phi^1_2 }^2
	+ 2 \norm{ \xi \phi^1_2 }^2 
	+ 2 \norm{ \overline e_1 e_1 \phi^1_2 }^2 
	+ 2 \norm{ e_1 \overline e_2 \phi^1_1 }^2 
	\nonumber \\
  & \quad
	+4 \norm{ e_1 \overline e_2 \phi^1_2 }^2 
	-4 \Real ( \overline e_2 e_2 \phi^1_2, 
		  \overline e_1 e_1 \phi^1_2 ) 
	+4 \Real ( i \xi \phi^1_2, e_1 \overline e_2 \phi^1_1 )
	\nonumber \\
  & \quad
	+ 2 \norm{ e_2 \overline e_1 \phi^2_2 }^2 
	+ 2 \norm{ \overline e_2 e_2 \phi^1_2 }^2
	+ 2 \norm{ \xi \phi^1_2 }^2  
	\nonumber \\
  & \quad
	+4 \Real (e_2 \overline e_1 \phi^2_2, i \xi \phi^1_2 )
	+ 4 \norm{ e_2 \overline e_1 \phi^1_2 }^2 
	- 4 \norm{ e_2 e_1 \phi^1_2 }^2
	\nonumber \\
  & \quad
	+ \norm{ \overline e_1 \overline e_1 \phi^2_2 }^2 
	+ \norm{ \overline e_2 \overline e_2 \phi^1_1 }^2  
	\nonumber \\
%
%
  &\sim  
	\norm{ e_1 e_1\phi^1_1 }^2 + \norm{ e_2 e_2\phi^2_2 }^2 
	+ 4 \norm{ e_1 \overline e_2 \phi^1_2 }^2 
	\nonumber \\
  &\quad 
	+ 4 \norm{ e_2 \overline e_1 \phi^1_2 }^2
	+ \norm{ \overline e_1 \overline e_1 \phi^2_2 }^2 
	+ \norm{ \overline e_2 \overline e_2 \phi^1_1 }^2  
	\nonumber \\
  &\quad
	+ \left( 
		2 \norm{ \overline e_2 e_2 \phi^1_2 }^2 
		- 4 \Real ( \overline e_2 e_2 \phi^1_2, 
		  	    \overline e_1 e_1 \phi^1_2 ) 
		+ 2 \norm{ \overline e_1 e_1 \phi^1_2 }^2 
	  \right)
	\label{eqn:key-penultimate}
	\\
  &\quad
	+ \left( 
		2 \norm{ \xi \phi^1_2 }^2 
		+ 4 \Real ( i \xi \phi^1_2, 
			    e_1 \overline e_2 \phi^1_1 )
		+ 2 \norm{ e_1 \overline e_2 \phi^1_1 }^2 
	  \right)
	\nonumber \\
  & \quad
	+ \left(
		2 \norm{ e_2 \overline e_1 \phi^2_2 }^2 
		+ 4 \Real (e_2 \overline e_1 \phi^2_2, 
			   i \xi \phi^1_2 )
		+ 2 \norm{ \xi \phi^1_2 }^2  
	  \right).
	\nonumber 
\end{align}
Now the three parts grouped in parentheses can be removed by the
Schwarz inequality.  This gives us 
\begin{align*}
\norm{ D^* \phi }^2 + 2 \norm{ \overline \partial_1 \phi }^2_1
  &\gtrsim  
	\norm{ e_1 e_1\phi^1_1 }^2 + \norm{ e_2 e_2\phi^2_2 }^2 
	+ 4 \norm{ e_1 \overline e_2 \phi^1_2 }^2 
	\\
  & \quad
	+ 4 \norm{ e_2 \overline e_1 \phi^1_2 }^2
	+ \norm{ \overline e_1 \overline e_1 \phi^2_2 }^2 
	+ \norm{ \overline e_2 \overline e_2 \phi^1_1 }^2  
\end{align*}
which is equation~\eqref{eqn:key-estimate}.
This concludes the proof of the Key Estimate. 
\end{proof}

Now to prove Theorem~\ref{thm:main-estimate}, we need an estimate
\begin{equation}
\label{eqn:subelliptic-estimate}
\norm{D^{\ast}\phi}^2 
	+ \norm{\overline \partial_1 \phi }_1^2
   \gtrsim c \norm{\phi}_2^2.
\end{equation}
In our local frame, the right-hand side of this equation can be
written as 
\begin{equation*}
c \norm{ \phi }_2^2 \sim
c \sum_{\alpha, \beta, j, k} 
	\norm{ e_{j} e_{k} \phi^{\alpha}_{\beta} }^2,
\end{equation*}
where $\alpha$ and $\beta$ run from $1$ to $2$, and $j, k \in \{
1, 2, \overline 1, \overline 2\}$.  We construct each of these
estimates individually, and organize them in the following lemma.
\begin{lemma}
\label{lemma:phi-est}
There exists a positive constant 
$C$ such that 
\begin{equation}
\label{eqn:phi-est}
\norm{D^{\ast}\phi}^2 
	+ \norm{\overline \partial_1 \phi }_1^2
	   \gtrsim C \norm{ e_{j} e_{k} \phi^{\alpha}_{\beta} }^2
\end{equation}
for all $j, k \in \{1,2,\overline 1,\overline 2\}$, $\alpha$, $\beta$, 
and $\phi \in \Gamma(M,E_2)$.  
\end{lemma}

\begin{proof}
What we will show, in fact, 
is that for each $j,k$ and each $\epsilon>0$ there is a constant
$C>0$ such that 
\begin{equation}
\label{eqn:phi-est-2}
\norm{D^{\ast}\phi}^2 
	+ \norm{\overline \partial_1 \phi }_1^2
    + \epsilon \norm{ \phi }_2^2
   \gtrsim C \norm{ e_{j} e_{k} \phi^{\alpha}_{\beta} }^2 .
\end{equation}
The constant $\epsilon$ can be chosen to be dominated by all the
different constants $C$, so that the sum of the various individual
estimates \eqref{eqn:phi-est} and \eqref{eqn:phi-est-2} yields the
subelliptic estimate~\eqref{eqn:subelliptic-estimate}. 

We prove this lemma in stages: we produce the
estimate~\eqref{eqn:phi-est} for each of the components $\phi^1_2$,
$\phi^2_1$, $\phi^1_1$, and $\phi^2_2$ in turn.  

{\sc The $\phi^1_2$ case:}
We begin by noting that we have the estimate~\eqref{eqn:phi-est} for
$\norm{ \overline e_2 e_1 \phi^1_2 }^2
	\sim \norm{ e_1 \overline e_2 \phi^1_2}^2 $ 
and 
$\norm{ \overline e_1 e_2 \phi^1_2 }^2
	\sim \norm{ e_2 \overline e_1 \phi^1_2 }^2 $ 
by the Key Estimate, Lemma~\ref{lemma:key-estimate}.

Now consider the part of inequality~\eqref{eqn:key-penultimate} that
we discarded in the last step of the proof of the Key Estimate:
\begin{eqnarray}
\lefteqn{ 
	\norm{ D^* \phi }^2 + 2 \norm{ \overline \partial_1 \phi }^2_1
} \nonumber \\
  & \gtrsim & 
	+ \left( 
		2 \norm{ \xi \phi^1_2 }^2 
		+ 4 \Real ( i \xi \phi^1_2, 
			    e_1 \overline e_2 \phi^1_1 )
		+ 2 \norm{ e_1 \overline e_2 \phi^1_1 }^2 
	  \right)
	\label{eqn:phi12-inequality}
	\\
  & & 
	+ \left(
		2 \norm{ e_2 \overline e_1 \phi^2_2 }^2 
		+ 4 \Real (e_2 \overline e_1 \phi^2_2, 
			   i \xi \phi^1_2 )
		+ 2 \norm{ \xi \phi^1_2 }^2  
	  \right).
	\nonumber 
\end{eqnarray}
Notice that, since $[\xi, e_j] \sim 0$ and $(i\xi)^* \sim i\xi$,  
\begin{align*}
\left|
	+ 4 \Real ( i \xi \phi^1_2, 
		    e_1 \overline e_2 \phi^1_1 )
\right|
& \sim 
	\left|
	+ 4 \Real ( e_2 \overline e_1 \phi^1_2, 
		    i \xi \phi^1_1 )
	\right| \\
& \lesssim 
	2 \left( 
		\frac{1}{\epsilon} \norm{ e_2 \overline e_1 \phi^1_2 }^2 
		+ \epsilon \norm{ i \xi \phi_1^1 }^2
	  \right) \\
& \lesssim  
	2 \left( 
		\frac{1}{\epsilon} \norm{ e_2 \overline e_1 \phi^1_2 }^2 
		+ \epsilon \norm{ \phi }^2_2
	  \right) 
\end{align*}
for any $\epsilon > 0$.  As we have already estimated 
$\norm{e_2 \overline e_1 \phi^1_2}^2$, this allows us to obtain an
estimate 
\begin{equation*}
c \Real (i\xi \phi^1_2, e_1 \overline e_2 \phi^1_1 )
	\lesssim \norm{ D^* \phi }^2 + \norm{ \overline \partial_1 \phi }^2_1
		 + \epsilon \norm{ \phi }_2^2
\end{equation*}
for some $c>0$.  Similarly, we can obtain an estimate 
\begin{equation*}
c \Real (i\xi \phi^1_2, e_2 \overline e_1 \phi^2_2 )
	\lesssim \norm{ D^* \phi }^2 + \norm{ \overline \partial_1 \phi }^2_1
		 + \epsilon \norm{ \phi }_2^2
\end{equation*}
for some $c>0$.   From these estimates and
inequality~\eqref{eqn:phi12-inequality}, we obtain estimates for 
$\norm{ i\xi \phi^1_2 }^2 $, 
$\norm{ e_1 \overline e_2 \phi^1_1 }^2 
	\sim \norm{ \overline e_2 e_1 \phi^1_1 }^2$,
and 
$\norm{ e_2 \overline e_1 \phi^2_2 }^2 
	\sim \norm{ \overline e_1 e_2 \phi^2_2 }^2$.

We again return to a term we discarded at the end of the proof of the
Key Estimate: we have
\begin{equation*}
\norm{ D^* \phi }^2 + 2 \norm{ \overline \partial_1 \phi }^2_1
 \gtrsim 
		2 \norm{ \overline e_2 e_2 \phi^1_2 }^2 
		- 4 \Real ( \overline e_2 e_2 \phi^1_2, 
		  	    \overline e_1 e_1 \phi^1_2 ) 
		+ 2 \norm{ \overline e_1 e_1 \phi^1_2 }^2.
\end{equation*}
We may rewrite part of this as 
\begin{align*}
|- 4 \Real ( \overline e_2 e_2 \phi^1_2, 
  	    \overline e_1 e_1 \phi^1_2 ) |
 & \sim  
	\left|
	- 4 \Real ( i \xi \phi^1_2, \overline e_1 e_1 \phi^1_2 )
	- 4 \Real ( e_2 \overline e_2 \phi^1_2, 
			\overline e_1 e_1 \phi^1_2) 
	\right| \\
 & \sim  
	 \left|
	    - 4 \Real ( i \xi \phi^1_2, \overline e_1 e_1 \phi^1_2 )
	    - 4 \norm{ e_1 \overline e_2 \phi^1_2 }^2 
	 \right| \\
 & \lesssim  
	 \left|
	     4 \Real ( i \xi \phi^1_2, \overline e_1 e_1 \phi^1_2 )
	 \right|
	 + 4 \norm{ e_1 \overline e_2 \phi^1_2 }^2 \\
\end{align*}
In the same way as above, we can control the inner product on the
right.  Since we have an estimate already for $\norm{ e_1 \overline e_2
\phi^1_2 }^2$, we get estimates for 
$\norm{ \overline e_2 e_2 \phi^1_2}^2$ and 
$\norm{ \overline e_1 e_1 \phi^1_2 }^2$.  

We can integrate by parts to write 
\begin{align*}
\norm{ e_1 e_2 \phi^1_2 }^2 
 & \sim  
	\norm{ e_1 \overline e_2 \phi^1_2 }^2
	- \Real ( i \xi \phi^1_2, \overline e_1 e_1 \phi^1_2 ) \\
 & \lesssim  
	\norm{ e_1 \overline e_2 \phi^1_2 }^2
	+ \norm{ i \xi \phi^1_2 }^2 
	+ \norm{ \overline e_1 e_1 \phi^1_2 }^2.
\end{align*}
The previous estimates for the terms on the the right-hand side of this
inequality then establish estimates for 
$\norm{ e_1 e_2 \phi^1_2 }^2  \sim \norm{ e_2 e_1 \phi^1_2 }^2$.   

We use the fact that 
$e_{\alpha} \overline e_{\alpha} \sim 
 \overline e_{\alpha} e_{\alpha} - i \xi$
to get 
$$
\norm{ e_{\alpha} \overline e_{\alpha} \phi^1_2 }^2 
 \lesssim
 2 \left( 
 	\norm{ \overline e_{\alpha} e_{\alpha}  \phi^1_2 }^2
 	+ \norm{ i \xi \phi^1_2 }^2
 \right),
$$
which gives us estimates for 
$\norm{ e_1 \overline e_1 \phi^1_2 }^2$
and $\norm{ e_2 \overline e_2 \phi^1_2 }^2$.

Using integration by parts, we get an equality 
\begin{equation*}
\norm{ \overline e_1 e_2 \phi^1_2 }^2 + \norm{ e_2 e_2 \phi^1_2 }^2 
\sim \norm{ e_1 e_2 \phi^1_2 }^2 + \norm{ \overline e_2 e_2 \phi^1_2 }^2.
\end{equation*}
We thus obtain an estimate on $\norm{ e_2 e_2 \phi^1_2 }^2$ from the
estimates on $\norm{ e_1 e_2 \phi^1_2 }^2$ and 
$\norm{ \overline e_2 e_2 \phi^1_2 }^2$.  Using this same
trick, we have 
\begin{equation*}
\norm{ e_1 e_1 \phi^1_2 }^2 + \norm{ \overline e_2 e_1 \phi^1_2 }^2 
\sim \norm{ \overline e_1 e_1 \phi^1_2 }^2 + \norm{ e_2 e_1 \phi^1_2 }^2,
\end{equation*}
and we get an estimate on $\norm{ e_1 e_1 \phi^1_2 }^2$.

Using Lemma~\ref{lemma:phi-lemma} for the local expression of $\overline
\partial_1 \phi$, we get 
\begin{align*}
\norm{ \overline \partial_1 \phi }_1^2 
  & \gtrsim 
	\norm{ \overline e_1 ( \overline e_1 \phi^1_2 
			   - \overline e_2 \phi^1_1) }^2 \\
  & \sim    
	\norm{ \overline e_1 \overline e_1 \phi^1_2 }^2 
	- 2 \Real (\overline e_1 \overline e_1 \phi^1_2,
		   \overline e_1 \overline e_2 \phi^1_1 ) 
	+ \norm{ \overline e_1 \overline e_2 \phi^1_1 }^2.
\end{align*}
On the other hand, 
\begin{align*}
| - 2 \Real (\overline e_1 \overline e_1 \phi^1_2,
	     \overline e_1 \overline e_2 \phi^1_1 ) |
 & \sim  
	| - 2 \Real (\overline e_1 e_2 \phi^1_2,
		     e_1 \overline e_1 \phi^1_1 ) | \\
 & \lesssim  
	\epsilon \norm{ \phi }^2_2 
	+ \frac{1}{\epsilon} \norm{ \overline e_1 e_2 \phi^1_2 }^2.
\end{align*}
Since we've already estimated $\norm{ \overline e_1 e_2 \phi^1_2 }^2$,
this gives us an estimate on 
$\norm{ \overline e_1 \overline e_1 \phi^1_2 }^2$
and  
$\norm{ \overline e_1 \overline e_2 \phi^1_1 }^2 \sim
\norm{ \overline e_2 \overline e_1 \phi^1_1 }^2$.
Similarly, we may use 
\begin{align*}
\norm{ \overline \partial_1 \phi }_1^2 
  & \gtrsim 
	\norm{ \overline e_2 ( \overline e_1 \phi^1_2 
			   - \overline e_2 \phi^1_1) }^2 
\end{align*}
and 
\begin{align*}
| - 2 \Real (\overline e_2 \overline e_1 \phi^1_2,
	     \overline e_2 \overline e_2 \phi^1_1 ) |
 & \lesssim  
	\epsilon \norm{ \phi }^2_2 
	+ \frac{1}{\epsilon} \norm{ e_1 \overline e_2 \phi^1_2 }^2
\end{align*}
to obtain estimates on 
$\norm{ \overline e_1 \overline e_2 \phi^1_2 }^2 \sim 
\norm{ \overline e_2 \overline e_1 \phi^1_2 }^2$.
This completes the proof of the $\phi^1_2$ case of
Lemma~\ref{lemma:phi-est}.

{\sc The $\phi^2_1$ case:}
Recall that $\phi^2_1 = \phi^1_2$ by Lemma~\ref{lemma:phi-lemma}, so this
case follows from the $\phi^1_2$ case.

{\sc The $\phi^1_1$ case:}
We begin by recalling that we have our estimate for 
$\norm{ e_1 e_1 \phi^1_1 }^2$
and 
$\norm{ \overline e_2 \overline e_2 \phi^1_1 }^2$ 
by the Key Estimate, Lemma~\ref{lemma:key-estimate}.  We also remark
that we have estimated 
$\norm{ \overline e_1 \overline e_2 \phi^1_1 }^2 \sim
\norm{ \overline e_2 \overline e_1 \phi^1_1 }^2$
in the proof the $\phi^1_2$ case of Lemma~\ref{lemma:phi-est}.  

In the proof of the Key Estimate, Lemma~\ref{lemma:key-estimate}, we
did not use the fact that 
\begin{equation*}
\norm{ \overline \partial_1 \phi }_1^2 
\gtrsim
	\norm{ e_1 (\overline e_1 \phi^1_2 - \overline e_2 \phi^1_1) }^2
	+ \norm{ e_2 (\overline e_1 \phi^1_2 - \overline e_2 \phi^1_1)
		}^2.
\end{equation*}
From this fact we have that 
\begin{align*}
\norm{ D^* \phi }^2 + \norm{ \overline \partial_1 \phi }_1^2 
   \gtrsim & 
	\norm{ e_1 \overline e_1 \phi^1_2 }^2 
	- 2 \Real (e_1 \overline e_1 \phi^1_2,
		   e_1 \overline e_2 \phi^1_1 ) 
 \\  & 
	+ \norm{ e_1 \overline e_2 \phi^1_1 }^2 
	+ \norm{ e_2 \overline e_1 \phi^1_2 }^2 
 \\  & 
	- 2 \Real (e_2 \overline e_1 \phi^1_2,
		   e_2 \overline e_2 \phi^1_1 )
	+ \norm{ e_2 \overline e_2 \phi^1_1 }^2 
\end{align*}
Using the same method as in the proof of the $\phi^1_2$ case, and noting
that we have estimates for all of the $\phi^1_2$ terms, we obtain
estimates for 
$\norm{ e_1 \overline e_2 \phi^1_1 }^2$
and $\norm{ e_2 \overline e_2 \phi^1_1 }^2$.

Using our integration by parts trick, we see that 
\begin{equation*}
\norm{ \overline e_1 e_1 \phi^1_1 }^2
	+ \norm{ e_2 e_1 \phi^1_1 }^2
	= \norm{ e_1 e_1 \phi^1_1 }^2
	  + \norm{ \overline e_2 e_1 \phi^1_1 }^2.
\end{equation*}
We have estimates for both terms on the right-hand side, so this gives
us estimates for 
$\norm{ \overline e_1 e_1 \phi^1_1 }^2$
and $\norm{ e_2 e_1 \phi^1_1 }^2 \sim \norm{ e_1 e_2 \phi^1_1 }^2$.

Now we produce an estimate for $\norm{ i\xi \phi^1_1}^2$.  We can write 
$i \xi \sim [\overline e_{\alpha}, e_{\alpha}]$ for $\alpha =1,2$, so
integration by parts yields 
\begin{align*}
\norm{ i \xi \phi^1_1 }^2 
  \sim & 
	( \overline e_1 e_1 \phi^1_1 - e_1 \overline e_1 \phi^1_1,
	  \overline e_2 e_2 \phi^1_1 - e_2 \overline e_2 \phi^1_1 ) \\
  \sim & 
	( \overline e_1 e_1 \phi^1_1, \overline e_2 e_2 \phi^1_1 )
	- ( e_1 \overline e_1 \phi^1_1, e_2 \overline e_2 \phi^1_1 )
 \\  & 
	- ( \overline e_1 e_1 \phi^1_1, \overline e_2 e_2 \phi^1_1 )
	+ ( e_1 \overline e_1 \phi^1_1, e_2 \overline e_2 \phi^1_1 ) \\
  \sim & 
	\norm{ e_1 e_2 \phi^1_1 }^2 
	- \norm{ e_1 \overline e_2 \phi^1_1 }^2 
	- \norm{ \overline e_1 e_2 \phi^1_1 }^2 
	+ \norm{ \overline e_1 \overline e_2 \phi^1_1 }^2  \\
  \lesssim &
	\norm{ e_1 e_2 \phi^1_1 }^2 
	+ \norm{ \overline e_1 \overline e_2 \phi^1_1 }^2 .
\end{align*}
This gives us an estimate on $\norm{ i \xi \phi^1_1 }^2$.

Since 
$e_1 \overline e_1 \phi^1_1 \sim \overline e_1 e_1 \phi^1_1 - i \xi
\phi^1_1$, we get 
$\norm{ e_1 \overline e_1 \phi^1_1 }^2 
\lesssim 2 \left( \norm{ \overline e_1 e_1 \phi^1_1 }^2 +
		  \norm{ i \xi \phi^1_1 }^2 \right)$
and an estimate on $\norm{ e_1 \overline e_1 \phi^1_1 }^2$.  Similarly, 
$\norm{ \overline e_2 e_2 \phi^1_1 }^2 
\lesssim 2 \left( \norm{ e_2 \overline e_2 \phi^1_1 }^2 +
		  \norm{ i \xi \phi^1_1 }^2 \right)$
and we may estimate $\norm{ \overline e_2 e_2 \phi^1_1 }^2$.  

Finally, integration by parts gives us the equalities
\begin{align*}
\norm{ \overline e_1 \overline e_1 \phi^1_1 }^2 
	+ \norm{ e_2 \overline e_1 \phi^1_1 }^2 
	& \sim 
		\norm{ e_1 \overline e_1 \phi^1_1 }^2 
		+ \norm{ \overline e_2 \overline e_1 \phi^1_1 }^2 \\
\intertext{and}
\norm{ e_2 e_2 \phi^1_1 }^2 
	+ \norm{ \overline e_1 e_2 \phi^1_1 }^2 
	& \sim 
		\norm{ \overline e_2 e_2 \phi^1_1 }^2 
		+ \norm{ e_1 e_2 \phi^1_1 }^2,
\end{align*}
which allow us to estimate 
$\norm{ \overline e_1 \overline e_1 \phi^1_1 }^2$, 
$\norm{ e_2 \overline e_1 \phi^1_1 }^2 \sim \norm{ \overline e_1 e_2 \phi^1_1
}^2$, 
and $\norm{ e_2 e_2 \phi^1_1 }^2$. 
This is the last of the required $\phi^1_1$ estimates, and so this
completes the proof of the $\phi^1_1$ case of
Lemma~\ref{lemma:phi-est}.

{\sc The $\phi^2_2$ case:}
It is simplest to notice the symmetry between the $\phi^2_2$ case and
the $\phi^1_1$ case.  For example, the Key Estimate gives us an
estimate on 
$\norm{ e_1 e_1 \phi_1^1 }^2$
and
$\norm{ e_2 e_2 \phi_2^2 }^2$
as well as 
$\norm{\overline e_2 \overline e_2 \phi_1^1 }^2$
and
$\norm{\overline e_1 \overline e_1 \phi_2^2 }^2$.
Making the appropriate changes in the proof of the $\phi^1_1$ case
will then give us a proof in this case as well.
As this is the final case, we have now completed the proof of
Lemma~\ref{lemma:phi-est}.   
\end{proof}

\section{A family of CR structures}
\label{section:versal-family}
In this section, we introduce an explicit family of CR structures 
parameterized by a finite-dimensional analytic set, and show that it
gives a local family of solutions to the deformation problem 
\begin{equation}
\label{eqn:deform-family}
\left.
\begin{split}
    P(\phi)		       & = 0, \\
    \overline\partial_0^* \phi & = 0.
\end{split}
\right\}
\end{equation}

We begin by saying precisely what we mean by a family of CR
structures.  Let $(M,{}^0T^{\prime\prime})$ be a compact strictly
pseudoconvex CR manifold of real dimension $2n-1$.   By a
\emph{family of deformations} of a given CR structure
${}^0T^{\prime\prime}$ we mean a triple
$(M,{}^{\phi(t)}T^{\prime\prime},T)$, where
$T \subset \C^{k}$ is a complex analytic subset containing the
origin $\origin$ and 
$\phi : T \to \Gamma(M,T'\otimes (^0T^{\prime\prime})^{\ast})$  
is a complex analytic map such that, for each $t\in T$,  $\phi(t)$
determines an integrable CR structure ${}^{\phi(t)}T^{\prime\prime}$ 
on $M$.  Recall that this means that $P(\phi(t)) = 0$ for all $t \in 
T$, as $P$ is the integrability condition for CR structures at finite 
distance from ${}^0T^{\prime\prime}$.  Finally, we require that 
$\phi(\origin) = 0$; that is, that $\phi(\origin)$ corresponds to the 
original CR structure ${}^0T^{\prime\prime}$. 
Then our main result of this section is the following theorem.
\begin{Theorem}
\label{thm:deform-existence}
Let $(M,{}^0T^{\prime\prime})$ be a compact, strictly pseudoconvex CR
manifold of real dimension $5$, and write
$\Harmonics = \ker \square$ for the set of harmonic elements of
$\Gamma(M,E_1)$.  Then there is a complex-analytic map 
$\phi\colon \Gamma(M,E_1) \to \Gamma(M,E_1)$ defined in a neighborhood
of $0$ such that if
\begin{equation}\label{eqn:def-T}
T=\{ t\in \Harmonics : 
	R_2(\phi(t)) = 
		\overline\partial_1 N\overline\partial_1^* L R_2 (\phi(t)) 
  \},
\end{equation}
then $(M,{}^{\phi(t)}T'',T)$ is a family of deformations of
${}^0T^{\prime\prime}$.  
\end{Theorem}

We will prove this theorem by constructing a locally complex analytic 
family of solutions to the deformation
problem~\eqref{eqn:deform-family}.  
We begin by producing some useful Sobolev estimates.

Our Laplacian $\square$ is a fourth-order differential operator, and so
we can expect that the Neumann operator gains four derivatives in the
directions of $\C \otimes H = {}^0T^{\prime} \oplus
{}^0T^{\prime\prime}$.  This is the content of the following lemma.
\begin{lemma} 
\label{lemma:Neumann-est}
Let $(M,{}^0T^{\prime\prime})$ be a compact, strictly pseudoconvex CR
manifold of dimension $5$.  For each integer $m \ge 0$, there exists a 
constant $c_m > 0$ such that  
\begin{equation*}
\norm{ N\psi }_{4,m} \leq c_m \norm{ \psi }_{0,m}
\end{equation*}
for all $\psi\in \Gamma(M,E_1)$.
\end{lemma}

\begin{proof}
We will show that 
\begin{equation}\label{eqn:4m-estimate}
\norm{ u }_{4,m} \leq c_m \norm{ \square u}_{0,m}
\end{equation}
whenever $u\in \Harmonics^\perp \cap \Gamma(M,E_1)$.
Because $\square$ is
subelliptic, $Nu$ is smooth whenever $u$ is smooth,
so the required estimate 
follows by approximating with smooth sections.

The proof of \eqref{eqn:4m-estimate} is
by induction on $m$.  By using a partition
of unity we may assume that $u$ is supported
in the domain of a frame satisfying \eqref{eqn:frame-comm}.
Observe that Lemma~\ref{lemma:Prop-8.1} and the Cauchy-Schwartz
inequality imply that
$\norm{u}\le \norm{\square u}$.
As usual, we will let $\sim$ and $\lesssim$ denote equality
and inequality modulo lower-weight terms, which can
be absorbed by using standard interpolation inequalities.  

We begin by considering derivatives in the $\xi$ direction.
By
Lemma~\ref{lemma:Prop-8.1} and Theorem~\ref{thm:main-estimate}, 
\begin{align*}
\norm{\xi u}_2^2
&\lesssim (\xi u, \square \xi u)\\
&\sim (\xi u, \xi \square u + [\square,\xi]u).
\end{align*}
Because $\xi$ commutes with $e_\alpha$ and $\overline e_\beta$ modulo
terms of weight $1$, it follows that $[\square,\xi]$ is an 
operator of weight at most $4$.  Therefore, after integrating by parts,
the second term above can be absorbed to 
 yield
\begin{equation}\label{eqn:xi-estimate}
\norm{\xi u}_2^2 \lesssim \norm{u}_4 \,\norm{\square u}
\lesssim  \epsilon \norm{u}_4^2+ \frac{1}{\epsilon}\norm{\square u}^2.
\end{equation}

Now we can prove \eqref{eqn:4m-estimate} for the case $m=0$.
Observe that the commutation relations for $e_\alpha$ and
$\overline e_\beta$ imply that $[e_\alpha , L]$ is equal to
a constant multiple
of $e_\alpha \xi$ modulo lower-weight terms.
Therefore,
using Lemma~\ref{lemma:Prop-8.1} and Theorem~\ref{thm:main-estimate}
again, we get 
\begin{align*}
\norm{u}_4^2 
&\sim \norm{Lu}_2^2\\
&\lesssim (Lu, \square Lu)\\
&\lesssim (Lu, L\square u) + (Lu, P_4 \xi u),
\end{align*}
where $P_4$ is some operator of weight $4$.
Integrating by parts and using \eqref{eqn:xi-estimate}, we find
\begin{equation*}
\norm{u}_4^2
\lesssim \norm{u}_4 \,\norm{\square u} + \norm{u}_4\,\norm{\xi u}_2,
\end{equation*}
so
\begin{align*}
\norm{u}_4
&\lesssim  \norm{\square u}
+ \norm{\xi u}_2\\
&\lesssim \epsilon\norm{u}_4 + \frac{1}{\epsilon} \norm{\square u}.
\end{align*}
Choosing $\epsilon$ small enough, 
we can absorb the $\norm{u}_4$ term and obtain 
\eqref{eqn:4m-estimate} when $m=0$.

Now assume that \eqref{eqn:4m-estimate} holds for some $m>0$.
By induction, we have
\begin{align*}
\norm{\xi u}_{4,m} 
&\lesssim \norm{\square \xi u}_{0,m}\\
&\lesssim \norm{\xi \square u}_{0,m} + \norm{[\square,\xi]u}_{0,m}\\
&\lesssim \norm{\square u}_{0,m+1} + \norm{u}_{4,m}\\
&\lesssim \norm{\square u}_{0,m+1}.
\end{align*}
If $e$ denotes any of the vector fields $e_\alpha$ or
$\overline e_\beta$, then $[\square,e] = P_3\xi + P_4$, where
$P_3$ and $P_4$ are operators of weight $3$ and $4$, respectively.
Thus
\begin{align*}
\norm{e u}_{4,m} 
&\lesssim \norm{\square e u}_{0,m}\\
&\lesssim \norm{e \square u}_{0,m} + \norm{[\square,e]u}_{0,m}\\
&\lesssim \norm{\square u}_{0,m+1} + \norm{\xi u}_{3,m}+ \norm{u}_{4,m} \\
&\lesssim \norm{\square u}_{0,m+1}.
\end{align*}
Since $\norm{u}_{4,m+1}$ is a sum of terms of the form and
$\norm{\xi u}_{4,m}$ and 
$\norm{e u}_{4,m}$,
this completes the induction.
\end{proof}

Recall that, for $\phi\in \Gamma(M,E_1)$, 
the almost CR structure ${}^{\phi}T^{\prime\prime}$ is integrable
exactly when $P(\phi)= \overline\partial_1\phi + R_2(\phi) = 0$.  
With this in mind, we state the following proposition (compare to 
\cite[Proposition 3.12, page 813]{Akahori:1978-family}). 
\begin{prop} 
\label{prop:R2-est}
Let $(M,{}^0T^{\prime\prime})$ be a compact, strictly pseudoconvex CR
manifold of dimension $5$.  Then for each positive integer $m \geq n$,
there exists a positive constant $\tilde{c}_m$ such that 
\begin{equation*}
\norm{ \overline\partial_1^{\ast} L R_2 (\phi) }_{0,m}
  \leq 
  \tilde{c}_m \norm{ \phi }_{4,m}^2
\end{equation*}
for all $\phi\in \Gamma(M,E_1)$.
\end{prop}
\begin{proof}
The proof of this proposition is simply the fact 
that $\overline \partial_1^*$, $L$, and $R_2$ take derivatives only in
the $\C \otimes H$
directions; thus 
$\overline\partial_1^{\ast} L R_2 (\phi)$ can be written in a local
frame for ${}^0T^{\prime}$ as a homogeneous quadratic polynomial in the
coefficients of $\phi$ and their derivatives, in which 
each monomial 
has a total of no more than four $\C\otimes H$ derivatives.  The
assumption that $m \ge n$ and the Sobolev embedding theorem yield
the result.
\end{proof}

Thus Proposition~\ref{prop:R2-est} combined with
Lemma~\ref{lemma:Neumann-est} in the case 
$\psi = \overline\partial_1^* L R_2 (\phi)$ 
yields the following theorem.
\begin{Theorem} 
\label{thm:NdR2-est}
Let $(M,{}^0T^{\prime\prime})$ be a compact, strictly pseudoconvex CR
manifold of dimension $5$.  
For each integer $m\ge n$, there exists a constant $\hat{c}_m > 0$
such that 
\begin{equation*}
\norm{ N \overline\partial_1^* L R_2 (\phi) }_{4,m}
 \leq \hat{c}_m \norm{ \phi }_{4,m}^2
\end{equation*}
for all $\phi\in \Gamma(M,E_1)$.
\end{Theorem}
We now use Theorem~\ref{thm:NdR2-est} to prove the main theorem of
this section, Theorem~\ref{thm:deform-existence}.  

\begin{proof}[Proof of Theorem~\ref{thm:deform-existence}]
We will solve this problem first in a Banach space: complete
$\Gamma(M,E_1)$ with respect to the norm $\norm{ \ }_{2,m}$ for some
integer $m\geq n$ to obtain a Banach space, which we denote by
$\Gamma_{2,m}(M,E_1)$.  Consider the Banach analytic map
from $\Gamma_{2,m}(M,E_1)$ to itself given by  
\begin{equation*}
\phi \mapsto 
	\phi + N \overline\partial^{\ast}_1 L R_2 (\phi).
\end{equation*}
Theorem~\ref{thm:NdR2-est} implies that $\phi \in \Gamma_{2,m}(M,E_1)$
is actually mapped to another element of $\Gamma_{2,m}(M,E_1)$.  
This is clearly an analytic local isomorphism.
The Banach inverse mapping theorem then gives us an analytic inverse
map; that is, an analytic function $s \mapsto \phi(s)$ from
$\Gamma_{2,m}(M,E_1)$ to itself such that 
\begin{equation}
\label{eqn:Banach-inverse-family}
\phi(s) + N \overline\partial^*_1 L R_2 (\phi(s)) = s, 
	\quad s\in \Gamma_{2,m}(M,E_1). 
\end{equation}

Our family \eqref{eqn:Banach-inverse-family} is locally (near the
origin $\origin$) parametrized by 
the analytic set $T$ defined in \eqref{eqn:def-T}.
To see this precisely, notice
that equation~\eqref{eqn:Banach-inverse-family} implies that 
for $t\in\Harmonics$,
\begin{equation}
\label{eqn:del-B-i-f}
\overline \partial_1 \phi(t) 
   + \overline \partial_1 N \overline \partial_1^* L R_2( \phi(t) )
   = 0
\end{equation}
(as $\overline \partial_1  = 0$ on $\Harmonics$).  Combining
this with the definition of $T$, we see that 
\begin{equation*}
T = \{ t\in \Harmonics : P(\phi(t)) = 0 \}.
\end{equation*}
Since $\phi(t)$ depends complex analytically on $t \in T$, our $T$ is
a complex analytic subset of $\Harmonics$.  
\end{proof}
\section{Proof of Versality}
\label{sect:versality}
In this section we prove that the family of CR structures constructed
in Theorem~\ref{thm:deform-existence} is versal, at least 
with respect to deformations of complex structure parametrized
by smooth complex manifolds.  In order to define
the notion of versality, we first make clear our definition of
deformations of a complex manifold $U$.  (In practice, $U$ will be a
complex neighborhood of our CR manifold $M$, which is embedded as a
hypersurface in a complex manifold $N$.)  A \emph{family of
deformations} of the complex manifold $U$ is a triple 
$(\mathcal{U},\pi,S)$, where 
$S \subset \C^{k}$ is a complex analytic subset containing the
origin $\origin$, $\mathcal{U}$ is a complex analytic
space that is differentiably (but not necessarily complex analytically)
isomorphic to $U\times S$, and 
$\pi : \mathcal{U} \equiv U\times S \to S$ is 
projection onto the second factor.  

We remark that a family of deformations $(\mathcal U,\pi ,S)$ of a
complex manifold $U$ gives rise to a unique $T'U$-valued 1-form
$\omega(s) \in \Gamma(U,T'U\otimes (T^{\prime\prime}U)^{\ast})$,
depending complex analytically on $S$.  Moreover, the complex
structure over $\pi ^{-1}(s)$, defined by 
\begin{displaymath}
{}^{\omega(s)}T^{\prime\prime}
  = \{ \overline X + \omega(s)(\overline X) \ : 
	\ \overline X \in T^{\prime\prime}U \},
\end{displaymath}
is integrable.  Conversely, by the Newlander-Nirenberg theorem, if
such an $\omega(s)$ is given, at least in the case in which
$S$ is nonsingular, then we can construct a family
of deformations $(\mathcal U,\pi,S)$ of the complex manifold $U$.  

Now suppose $(M,{}^0T^{\prime\prime})$ is a strictly
pseudoconvex CR manifold.  A family of deformations
$(M,^{\phi(t)}T^{\prime\prime},T)$ of CR structures over M is said to
be \emph{versal} if 
whenever $(M,{}^0T^{\prime\prime})$ is embedded as
a real hypersurface in an $n$-dimensional complex manifold $N$ and 
$(\mathcal U,\pi,S)$ is any deformation of the complex structure
on a neighborhood
$U$ of $M$ in $N$, we have the following two conditions.  First,
there exists a neighborhood of the origin $S'\subset S$ for which 
there is a holomorphic map $h : S'\rightarrow T$ 
and smooth embeddings $f(s) : M \rightarrow \pi^{-1}(s)$ for all 
$s \in S'$ such that $h(\origin) = \origin$ and $f(\origin)$ is the
identity map.  Second, we note that $\omega(s)$ induces a CR structure
over $M$ when we consider $M$ embedded in $U$ via $f(s)$. 
Let us denote this CR structure by ${}^{\omega(s)\cdot f(s)}T''$.
If $s$ is sufficiently close to the origin, this defines a unique
deformation tensor $\omega(s)\cdot f(s)\in 
\Gamma(M,T'\otimes (^0T^{\prime\prime})^{\ast})$ 
by 
\begin{equation}
\label{eqn:induced-CR}
{}^{\omega(s)\cdot f(s)}T''
 = \left\{
	\overline X + (\omega(s)\cdot f(s) ) (\overline X) 
	: \overline X \in {}^0T^{\prime\prime}
   \right\}.
\end{equation}
Our requirement is that
this CR structure be the same as the one induced by $\phi$ at the
point $h(s) \in T$:
\begin{displaymath}
\omega(s)\cdot f(s)=\phi(h(s)) \ \text{ for all $s\in S'$.}
\end{displaymath}
We will only deal with smooth deformations; that is, deformations in
which the analytic space $S$ is, in fact, a complex manifold, rather
than a variety with singularities.

We now state our main theorem of this section.  
\begin{Theorem} 
\label{thm:10-1}
Suppose $(M,{}^0T^{\prime\prime})$ is a compact strictly pseudoconvex CR
manifold of real dimension $2n-1 = 5$ that is embedded as a real
hypersurface in a complex manifold $N$ of complex dimension $n=3$.  If
the family of CR deformations $(M,{}^{\phi(t)}T^{\prime\prime},T)$ is
a smooth family of deformations, then it is versal with respect to
smooth deformations (that is, with respect to deformations
$(\mathcal{U},\pi,S)$ of a neighborhood $U$ of $M$ in $N$, where the
analytic space $S$ is a complex manifold). 
\end{Theorem}
Our proof can be modified to work in the case that $S$ has a
singularity, so the claim would be that the family of CR deformations
is versal.  We leave this claim to another paper.

\begin{proof}
We must construct $h(s)$ and $f(s)$.  Suppose that we are
given a family of deformations of a neighborhood $U$ of $M$,
$(\mathcal U,\pi,S)$.  Let $\{U_j\}$ be a covering of $U$ by 
coordinate domains,
indexed by some finite set.  Let $\{z_j^1,z_j^2,z_j^3\}$
be local holomorphic coordinates on $U_j$, and let 
$\tau_{jk}^{l}(z_k^1,z_k^2,z_k^3)$ be transition functions:
\begin{displaymath}
z_j^{l} = \tau_{jk}^{l}(z_k^1,z_k^2,z_k^3), 
 \ l= 1,2,3, \ \text{on} \ U_j\cap U_k. 
\end{displaymath}
For brevity, we will write this as
\begin{displaymath}
z_j = \tau_{jk}(z_k), \ \text{on} \ U_j\cap U_k.
\end{displaymath}
We can extend this to a local coordinate covering 
$\{ U_j\times S \}$ for $(\mathcal U,\pi,S)$ 
with transition functions 
$\tau_{jk}^{l}(z_k^1,z_k^2,z_k^3,s)$ 
defined on $U_j\times S \cap U_k\times S$, holomorphic in $z^j_k$ and
smooth in $s$.  
We use a similar abbreviation as above:
\begin{displaymath}
z_j = \tau_{jk}(z_k,s), 
\ \text{on} \ U_j\times S\cap U_k\times S,
\end{displaymath}
with the requirement that  $\tau_{jk}(z_k,\origin)=\tau_{jk}(z_k)$.  
For simplicity, we use local complex coordinates 
$\{ z_j^1(s),z_j^2(s),z_j^3(s) \}$ depending complex analytically
on the parameter $s$.  That is, each function $z_j^k(s)$ is a smooth
function on $U_j$ and complex analytic on $S$, and the corresponding
complex structure on $\pi^{-1}(s)$ (as an element of
$\Gamma(U,T'U\otimes (T^{\prime\prime}U)^{\ast})$ )
is determined by
\begin{displaymath}
(\overline X+\omega(s)(\overline X)) z_j^k(s)=0, 
  \ \text{for all} \ \overline X\in T^{\prime\prime}U.
\end{displaymath}
Similarly, the induced CR structure defined in
equation~\eqref{eqn:induced-CR} is also determined locally by  
\begin{equation}
\label{eqn:induced-CR-X}
(\overline X+ \omega(s)\cdot f_j(s) (\overline X)) f_j^{l}(s)=0 
  \ \text{for all} \ \overline X\in {}^{0}T^{\prime\prime},
\end{equation}
where $f_j^{l}(s) = z_j^{l} \circ f(s)$.  
This equality also means that the map $f(s)$ is a CR
embedding from $(M,{}^{\omega(s)\cdot f(s)}T^{\prime\prime})$ to
$\pi^{-1}(s)$, with the complex structure $\omega(s)$.

We have to construct $f(s)$, locally expressed by
$f_j(s)=(f_j^1(s),f_j^2(s),f_j^3(s))$ on $U_j$, which depends complex
analytically on $S$, and a holomorphic map $h$ from $S$ to 
$T \subset \Harmonics$, satisfying
\begin{eqnarray*}
f_j(s)                  & = & \tau_{jk}(f_k(s),s) \\
\omega(s) \cdot f_j(s)  & = & \phi(h(s)) 
\end{eqnarray*}
for all $s \in S$ (where, if necessary, we may shrink $S$ to a smaller
neighborhood of $\origin$).  
The proof of the existence of such functions is a standard formal
power series argument.  Consider the power series expansions
\begin{displaymath}
f_j(s) = \sum_{|\alpha| = 0}^{\infty} f_{j\mid \alpha} s^{\alpha}
\qquad \text{ and } \qquad
h(s)   = \sum_{|\alpha| = 0}^{\infty} h_{\mid \alpha}  s^{\alpha}.
\end{displaymath}
We are using multi-index notation, so if $s = (s_1,\ldots,s_r)$
and $\alpha = (\alpha_1,\ldots,\alpha_r)$, then $|\alpha| = \alpha_1 +
\cdots + \alpha_r$ and $s^{\alpha} = s_1^{\alpha_1} \cdots
s_r^{\alpha_r}$.  
In general, if $F$ is any vector-bundle-valued
function of $s$, we will use the notation 
$\kappa_{m}F$ to mean the part of the power series for $F(s)$
about $s=0$ that is homogeneous of order $m$ in $s$.   
For such homogeneous polynomials, we will use a subscript $(k)$ to
indicate the degree in $s$.  Similarly, a superscript $(k)$ will
indicate a (not usually homogeneous) polynomial of degree $k$ in $s$.

First we formally construct these power series, then prove
convergence.  Let $f_j^{(m)}$ and $h^{(m)}$ be the $m$th partial sums
in the above power series expansions:
\begin{displaymath}
f_j^{(m)}(s) = \sum_{|\alpha|=0}^{m} f_{j\mid \alpha} s^{\alpha}
\qquad \text{ and } \qquad
h^{(m)}(s)   = \sum_{|\alpha|=0}^{m} h_{\mid \alpha}  s^{\alpha}.
\end{displaymath}
We construct $f_j^{(m)}(s)$ and $h^{(m)}(s)$ formally by induction on $m$.

At any step $m$, we wish to have $f_j^{(m)}$ and $h^{(m)}$ satisfy 
\begin{equation}
\label{eqn:inductive-series}
\left.
\begin{split}
f_j^{(m)}(s) & = \tau_{jk} (f_k^{(m)}(s),s) +O(|s|^{m+1}) \\
\omega(s) \cdot f_j^{(m)}(s)
	     & = \phi(h^{(m)}(s)) +O(|s|^{m+1})
\end{split}
\right\}
\end{equation}
for $s \in S$ near $\origin$.

At our initial step (that is, at $m=0$), we define 
$f_j^{(0)}(s)= z_j(s)$ and $h^{(0)}(s)=0$.
These obviously satisfy our criterion~\eqref{eqn:inductive-series}.

Now we assume that we have already constructed $f_j^{(m)}$ and
$h^{(m)}$ satisfying~\eqref{eqn:inductive-series}.  To begin our
construction of $f_j^{(m+1)}$ and $h^{(m+1)}$, we define a polynomial
$g_{j\mid (m+1)}$ on $U_j$, homogeneous of degree $m+1$ in $s$, such
that 
\begin{equation}
\label{eqn:g-properties}
f_j^{(m)}(s) + g_{j\mid (m+1)}(s)
  = \tau_{jk}( f_k^{(m)}(s) + g_{k\mid (m+1)}(s) , s )
	 +O(|s|^{m+2}).
\end{equation}
(In this way, $g_{j\mid (m+1)}$ is a rough first approximation of
$\kappa_{m+1}(f^{(m+1)}_j)$, the homogeneous part of
$f^{(m+1)}_j$ in degree $m+1$.)  
To do this, we construct vector-valued polynomials $\sigma_{jk\mid (m+1)}$
on $U_j\cap U_k$, again homogeneous of degree $m+1$ in $s$, 
by the relation
\begin{equation}
\label{eqn:sigma-defn}
\sigma_{jk\mid (m+1)}(s)  
 =\tau_{jk}( f_k^{(m)}(s),s ) -  f_j^{(m)}(s) 
	+O(|s|^{m+2}).
\end{equation}
This definition of $\sigma_{jk\mid (m+1)}$ makes sense as the
induction hypothesis~\eqref{eqn:inductive-series} 
implies that the right-hand side of equation~\eqref{eqn:sigma-defn}
has only terms of order $m+1$ and higher in $s$.
We use these $\sigma_{jk\mid (m+1)}$ and a partition of unity 
$\{ \rho_j \}$ subordinate to the covering $\{ U_j \}$ to define 
\begin{eqnarray}\label{eqn:def-g}
g_{j \mid (m+1)}(s) = \sum_{k} \rho_k \sigma_{jk\mid (m+1)}(s).
\end{eqnarray}
We will show that such $g_{j \mid (m+1)}$
satisfy~\eqref{eqn:g-properties}.
To do this, we need to know how $g_{j \mid (m+1)}$ (or $\sigma_{jk 
\mid (m+1)}$) transform over different coordinate charts.  We have the
following lemma (compare to 
\cite[Lemma 3.2, page 828]{Akahori-Miyajima:1980}).

\begin{lemma} \label{lemma:sigma-trans}
On $U_j\cap U_k\cap U_{l}$,
\begin{equation}
\label{eqn:sigma-trans}
\sigma_{jk \mid (m+1)}(s)
  + \frac{\partial \tau_{jk}} {\partial z_{k}}(f_k^{(m)}(s),s)
    \sigma_{kl \mid (m+1)}(s)
  = \sigma_{jl \mid (m+1)}(s) 
	+O(|s|^{m+2}).
\end{equation}
\end{lemma}

\begin{proof}
By the definition of $\sigma_{jk\mid (m+1)}$,
\begin{displaymath}
\sigma_{jk\mid (m+1)}(s) 
	= \tau_{jk}(f_k^{(m)}(s),s)
-f_j^{(m)}(s) 	+O(|s|^{m+2}).
\end{displaymath}
We replace $f_k^{(m)}(s)$ 
with 
$\tau_{kl}(f_{l}^{(m)}(s),s)- \sigma_{kl \mid (m+1)}(s)$
to get 
\begin{displaymath}
\sigma_{jk\mid (m+1)}(s) 
	= \tau_{jk}(\tau_{kl}(f_{l}^{(m)}(s),s)
			   - \sigma_{kl \mid (m+1)}(s) 
			   ,s)
-
f_j^{(m)}(s) 
	        	+O(|s|^{m+2}).
\end{displaymath}
We expand the first term on the right-hand side in a power series
about the point $(z_k,s) = (\tau_{kl}(f_{l}^{(m)}(s),s),s)$; this
implies  
\begin{equation}
\label{eqn:sigma-trans-proof}
\begin{aligned}
\sigma_{jk\mid (m+1)}(s) 
	 & = \tau_{jk}(\tau_{kl}(f_{l}^{(m)}(s),s),s) 
	        - \frac{\partial \tau_{jk}}{\partial z_k}
	   (\tau_{kl}(f_{l}^{(m)}(s),s),s)\sigma_{kl \mid (m+1)}(s) \\
&\quad-f_j^{(m)}(s) 
	+O(|s|^{m+2})\\
& = \tau_{jl}(f_{l}^{(m)}(s),s)
	        - \frac{\partial \tau_{jk}}{\partial z_k}
	   (f_k^{(m)}(s),s)\sigma_{kl \mid (m+1)}(s) \\
&\quad-f_j^{(m)}(s) 
	+O(|s|^{m+2}).\\
\end{aligned}
\end{equation}
(In the last line we have used the inductive hypothesis
\eqref{eqn:inductive-series} and Taylor's theorem
applied
to $\partial \tau_{jk}/\partial z_k$.
Any error term involving $\sigma_{kl\mid(m+1)}(s)$
multiplied by itself or by $O(|s|^{m+1})$ can be
absorbed into
$O(|s|^{m+2})$.)
The first and third terms simplify
to $\sigma_{jl \mid (m+1)}(s)$ modulo $O(|s|^{m+2})$, and so
equation~\eqref{eqn:sigma-trans-proof} reduces to
equation~\eqref{eqn:sigma-trans}.  This proves the lemma.
\end{proof}

\begin{lemma}
With $g_{j\mid (m+1)}$ defined by \eqref{eqn:def-g}, 
$f_j^{(m)} + g_{j\mid (m+1)}$ transforms as in
equation~\eqref{eqn:g-properties}.  
\end{lemma}

\begin{proof}
From the definition of $g_{j\mid (m+1)}$ and 
\eqref{eqn:sigma-trans},
\begin{equation}\label{eqn:expand-g}
\begin{aligned}
g_{j\mid (m+1)} (s) 
&= \sum_l \rho_l \sigma_{jl\mid(m+1)}(s)\\
&= \sigma_{jk\mid(m+1)}(s) 
+ \sum_l \rho_l \frac{\partial \tau_{jk}}{\partial z_k} 
(f_k^{(m)}(s),s)\sigma_{kl\mid(m+1)}(s) + O(|s|^{m+2})\\
&= \sigma_{jk\mid(m+1)}(s)
+ \frac{\partial \tau_{jk}}{\partial z_k} 
(f_k^{(m)}(s),s)g_{k\mid (m+1)}(s)+ O(|s|^{m+2}).
\end{aligned}
\end{equation}
Thus
\begin{equation*}
\begin{aligned}
f_j&^{(m)}(s) + g_{j\mid (m+1)} (s) \\
&\quad= f_j^{(m)}(s)+
\sigma_{jk\mid(m+1)}(s)
+ \frac{\partial \tau_{jk}}{\partial z_k} 
(f_k^{(m)}(s),s)g_{k\mid (m+1)}(s)+ O(|s|^{m+2})\\
&\quad= \tau_{jk}(f_k^{(m)}(s),s) 
+ \frac{\partial \tau_{jk}}{\partial z_k} 
(f_k^{(m)}(s),s)g_{k\mid (m+1)}(s)+ O(|s|^{m+2}).
\end{aligned}
\end{equation*}
By Taylor's theorem, this is equivalent to \eqref{eqn:g-properties}.
\end{proof}

To define the next term in our formal power series, we will
write locally 
\begin{equation} 
\label{eqn:f-h-defn}
\begin{aligned}
f_j^{(m+1)}(s) & = f_j^{(m)}(s) + g_{j\mid(m+1)}(s) + \zeta_{j\mid(m+1)}(s) \\ 
h^{(m+1)}(s) & = h^{(m)}(s) + h_{(m+1)}(s),
\end{aligned}
\end{equation}
where $\zeta_{j\mid(m+1)}$
is the local expression for a 
homogeneous polynomial
$\zeta_{(m+1)}$ 
of degree $m+1$ in 
$s$, with values in 
$\Gamma(M,T')$, and
$h_{(m+1)}$ is  a homogeneous polynomial of degree $m+1$ with 
values in $\Harmonics$.
Since the transformation law for sections of $T'$ is
\begin{displaymath}
\zeta_{j \mid (m+1)}
  = \frac {\partial \tau_{jk}} {\partial z_k} \zeta_{k \mid (m+1)},
\end{displaymath}
it follows that 
our prospective $f^{(m+1)}(s)$ transforms the correct way:
\begin{multline*}
f_j^{(m)}(s) + g_{j\mid (m+1)}(s) + \zeta_{j \mid (m+1)}(s)\\
    = \tau_{jk} ( f_k^{(m)}(s) 
		       + g_{k\mid (m+1)}(s) 
		       + \zeta_{k \mid (m+1)}(s)
 		      , s) 
    +O(|s|^{m+2}).
\end{multline*}

We still must construct $\zeta_{(m+1)}$ and $h_{(m+1)}$ so that 
$f^{(m+1)}_j(s)$ and $h^{(m+1)}(s)$, defined as in
equation~\eqref{eqn:f-h-defn}, satisfy the inductive
hypothesis~\eqref{eqn:inductive-series}.  Note first that, by
equation~\eqref{eqn:induced-CR-X}, the CR structure defined by
$f^{(m+1)}(s)$ must satisfy 
\begin{align*}
(\overline X 
  + & (\omega(s) \cdot (f_j^{(m)}(s) 
		+ g_{j\mid (m+1)}(s)
		+ \zeta_{j\mid (m+1)}(s) 
    		)(\overline X)) \nonumber \\
    & (f_j^{(m)}(s) + g_{j\mid (m+1)}(s) + \zeta_{j\mid (m+1)}(s) ) = 0.
\end{align*}
From this it follows that 
\begin{align*}
\omega(s) \cdot & (f_j^{(m)}(s) 
		   + g_{j\mid (m+1)}(s) 
		   + \zeta_{j\mid (m+1)}(s) ) \\
& = \omega(s) \cdot (f_j^{(m)}(s) + g_{j\mid (m+1)}(s) ) 
	 + \overline\partial_{T'} \zeta_{j \mid (m+1)}(s) 
	 +O(|s|^{m+2}).
\end{align*}
On the other hand, from the definition, it is clear that (see
equation~\eqref{eqn:Banach-inverse-family}) the map $\phi$ linearizes
to the identity, so  
\begin{displaymath}
\phi( h^{(m)}(s) + h_{(m+1)}(s) ) 
 = \phi( h^{(m)}(s) ) + h_{(m+1)}(s) 
 +O(|s|^{m+2}).
\end{displaymath}
Finding solutions to the second equation in
\eqref{eqn:inductive-series} is thus reduced to the following theorem.
\begin{Theorem} \label{thm:10-2}
There are vector-valued polynomials $\zeta_{(m+1)}$ and $h_{(m+1)}$,
homogeneous in $s$ of degree $m+1$, solving 
\begin{equation}
\label{eqn:10-21}
\begin{aligned}
\omega(s) \cdot & ( f_j^{(m)}(s) + g_{j\mid (m+1)}(s) ) 
	+ \overline\partial_{T'} \zeta_{(m+1)}(s)  \\
 &  = \phi(h^{(m)}(s)) + h_{(m+1)}(s) +O(|s|^{m+2}),
\end{aligned}
\end{equation}
where $\zeta_{(m+1)}$ takes values in $\Gamma(M,T')$, and $h_{(m+1)}$
takes values in the finite-dimensional harmonic space $\Harmonics
\subset \Gamma(M,E_1)$. 
\end{Theorem}

The proof of this theorem will follow from several lemmas and 
propositions.  

\begin{prop} 
\label{prop:10-1}
There is a homogeneous polynomial 
$\theta_{(m+1)}$ of degree $m+1$ in $s$, 
with values in 
$\Gamma(M,{}^0T^{\prime})$,
such that  
\begin{displaymath}
\omega(s) \cdot (f_j^{(m)}(s)
  + g_{j\mid (m+1)}(s) 
  + \theta_{j\mid (m+1)}(s) )
  \in \Gamma(M,E_1),
\end{displaymath}
where we have written $\theta_{j\mid (m+1)}$ for 
$\theta_{(m+1)}\mid_{U_j}$. 
\end{prop}
\begin{proof}
Since our CR structure is strictly pseudoconvex, the map 
\begin{align*}
\Gamma(M,{}^0T^{\prime}) 
	& \rightarrow \Gamma(M,F\otimes (^0T^{\prime\prime})^{\ast}) \\
u 	& \mapsto \pi_F \overline\partial_{T'} u
\end{align*}
is an isomorphism.  Hence there is a 
$\Gamma(M,^0T^{\prime})$-valued polynomial $\theta$ 
which such that 
$\kappa_{m+1}(
	\omega(s)\cdot  
		(f_j^{(m)}(s)+g_{j\mid (m+1)}(s))  
	+ \overline\partial_{T'} \theta (s))$
is a polynomial that takes values in 
$\Gamma(M,{}^0T^{\prime}\otimes (^0T^{\prime\prime})^{\ast})$.
By the inductive hypothesis, for each $l < m$, the polynomial 
$\kappa_{l}( \omega(s)\cdot (f_j^{(m)}(s)+g_{j\mid (m+1)}(s)) )
 = \kappa_{l}( \omega(s)\cdot f_j^{(m)}(s) )$
already takes values in 
$\Gamma(M,{}^0T^{\prime}\otimes (^0T^{\prime\prime})^{\ast})$;
thus we may assume $\theta = O(|s|^{m+1})$.
Writing $\theta_{(m+1)}$ for $\kappa_{m+1} \theta$ and $\theta_{j\mid
(m+1)}$ for $\theta_{(m+1)}\mid_{U_j}$, we thus have 
\begin{displaymath}
\kappa_{m+1} (\omega(s) \cdot (f_j^{(m)}(s) 
				+ g_{j\mid (m+1)} (s)
				+ \theta_{j\mid (m+1)} (s))
	     ) 
\in \Gamma(M,{}^0T'\otimes(^0T^{\prime\prime})^{\ast})
\end{displaymath}
To prove the proposition, it suffices to show
\begin{equation}
\label{eqn:prop-necc}
\kappa_{m+1}(
	\pi_{F}
	\overline\partial^{(1)}
	(\omega(s) \cdot (
		f_j^{(m)}(s) + g_{j\mid (m+1)}(s) + \theta_{j\mid (m+1)}(s) 
			 )
	)   )
	= 0.
\end{equation}
In order to show this, we first prove the next lemma.
\begin{lemma} 
\label{lemma:10-2}
\begin{equation}
R_2(\phi(h^{(m)}(s)))
  = R_2(\omega(s) \cdot (f_j^{(m)}(s)))
  +O(|s|^{m+2})
\end{equation}
and
\begin{eqnarray}
R_3(\phi(h^{(m)}(s))) 
  = R_3(\omega(s) \cdot (f_j^{(m)}(s)))
  +O(|s|^{m+2})
\end{eqnarray}
hold.  In particular, 
$R_3(\omega(s)\cdot (f_j^{(m)}(s))) = O(|s|^{m+2})$
as $\phi(t) \in \Gamma(M,E_1)$.  
\end{lemma}
\begin{proof}
For $\psi \in \Gamma(M,T'\otimes (^0T^{\prime\prime})^{\ast})$, 
$R_k(\psi)$ ($k=2,3$) are the parts of the deformation equation that
are order $k$ in $\psi$.  (Of course, each $R_k(\psi)$ includes first 
derivatives of $\psi$.)  The expressions for $R_k$ are given in
equations~\eqref{eqn:R2-defn} and~\eqref{eqn:R3-defn}.
Since $R_2$ is quadratic, and we may replace each $\phi(h^{(m)}(s))$
with $\omega(s) \cdot (f_j^{(m)}(s))$ in turn.  On the one hand,
$\phi(h^{(m)}(s)) = \omega(s) \cdot (f_j^{(m)}(s)) +O(|s|^{m+1})$ by the induction hypothesis~\eqref{eqn:inductive-series}.
On the other hand, $\phi(h^{(m)}(s))$ itself satisfies 
$\phi(h^{(m)}(s)) =O(|s|)$.  Together, these facts imply
that 
$R_2(\phi(h^{(m)}(s))) = R_2(\omega(s) \cdot (f_j^{(m)}(s)))
+O(|s|^{m+2})$.  The proof for $R_3$ is similar.
\end{proof}

Continuing the proof of Proposition~\ref{prop:10-1}, 
we remark that $\omega(s)$ is, for each $s$, an integrable complex
structure.   Since $(f_j^{(m)}(s) + g_{j\mid (m+1)} + \theta_{j\mid 
(m+1)})$ is a CR embedding for each $s$, modulo terms of order $m+2$
and higher, the CR structure induced by $\omega(s)$ is also integrable:
\begin{displaymath}
\left( \overline\partial^{(1)} + R_2 + R_3  \right) 
	\omega(s) \cdot 
		(f_j^{(m)}(s) + g_{j\mid (m+1)}(s)
		+ \theta_{j\mid (m+1)}(s) )
 =O(|s|^{m+2}).
\end{displaymath}
Obviously, we may remove the terms of order $m+2$ and higher to see
that 
\begin{multline*}
R_2(\omega(s) \cdot 
	(f_j^{(m)}(s) + g_{j\mid (m+1)}(s) + \theta_{j\mid (m+1)}(s) ) )\\
  =  R_2( \omega(s) \cdot (f_j^{(m)}(s) ) ) +O(|s|^{m+2}).
\end{multline*}
From the previous lemma, 
$R_2( \omega(s) \cdot (f_j^{(m)}(s) ) ) = R_2( \phi(h^{(m)}(s)) )
+O(|s|^{m+2})$, and so 
\begin{displaymath}
R_2(\omega(s) \cdot 
	(f_j^{(m)}(s) + g_{j\mid (m+1)}(s) + \theta_{j\mid (m+1)}(s) ) )
  =  R_2( \phi(h^{(m)}(s)) ) +O(|s|^{m+2}).
\end{displaymath}
A similar computation shows that 
\begin{displaymath}
R_3 (\omega(s) \cdot 
	(f_j^{(m)}(s) + g_{j\mid (m+1)}(s) + \theta_{j\mid (m+1)}(s) ) )
  =O(|s|^{m+2})
\end{displaymath}
(and the zero follows from Lemma~\ref{lemma:10-2}).  
The integrability condition is thus 
\begin{align*}
\overline\partial^{(1)}
	&\omega(s)\cdot 
		(f_j^{(m)}(s) 
		 + g_{j\mid (m+1)}(s)
		 + \theta_{j\mid (m+1)}(s))
	 + R_2(\phi(h^{(m)}(s))) \\
	& =O(|s|^{m+2}).
\end{align*}
Because $\phi(t)$ takes its values in $\Gamma(M,E_1)$, we have 
$\pi_{F} (R_2(\phi(h^{(m)}(s)))) = 0$.
Hence
\begin{displaymath}
\pi_F (\overline\partial^{(1)}
	(\omega(s)\cdot 
		(f_j^{(m)}(s) 
		 + g_{j\mid (m+1)}(s)
		 + \theta_{j\mid (m+1)}(s) 
		)
	)
      )
 	=O(|s|^{m+2}).
\end{displaymath}
This is equivalent to equation~\eqref{eqn:prop-necc}, so this proves
Proposition~\ref{prop:10-1}.    
\end{proof}

\begin{lemma} \label{lemma:10-3}
\begin{displaymath}
(1 - \overline\partial_1 N \overline\partial_1^{\ast} L )
	R_2(\omega(s)\cdot ( f_j^{(m)}(s) 
			       + g_{j\mid (m+1)}(s)
 			       + \theta_{j\mid (m+1)}(s)
			     )
	     )
	=O(|s|^{m+2}).
\end{displaymath}

\end{lemma}
\begin{proof}
We recall that 
\begin{displaymath}
P (\omega(s) \cdot 
	(f_j^{(m)}(s) 
	 + g_{j\mid (m+1)}(s)
	 + \theta_{j\mid (m+1)}(s)
	)
  ) =O(|s|^{m+2}).
\end{displaymath}
(The map defined on each $U_j$ by 
$f_j^{(m)}(s) + g_{j\mid (m+1)}(s) + \theta_{j\mid (m+1)}(s)$ makes sense
globally modulo $O(|s|^{m+2})$.)  Thus
\begin{align*}
\overline\partial_1 \omega(s) &\cdot 
	(f_j^{(m)}(s)+g_{j\mid (m+1)}(s)+\theta_{j\mid (m+1)}(s))) \\
	&+ R_2(\omega(s)\cdot (f_j^{(m)}(s) 
	+ g_{j\mid (m+1)}(s)+\theta_{j\mid (m+1)}(s))))   
	=O(|s|^{m+2}).
\end{align*}
We apply the operator
$1-\overline\partial_1 N \overline\partial_1^{\ast} L$ to this
equality.  By Proposition~\ref{prop:10-1}, the left-hand side is the
image of an element of $\Gamma(M,E_1)$ under $\overline \partial_1 +
R_2$, so this makes sense.  The decomposition of
Theorem~\ref{thm:Hodge-Kodaira} implies that 
$(1-\overline\partial_1 N \overline \partial_1^{\ast} L)
\overline\partial_1 =0$, and from this Lemma~\ref{lemma:10-3} follows
easily. 
\end{proof}

\begin{prop} 
\label{prop:10-2}
\begin{displaymath}
\overline\partial_1 
\left[
 \omega(s) \cdot 
	( f_j^{(m)}(s) + g_{j\mid (m+1)}(s) + \theta_{j \mid (m+1)}(s) ) 
 - \phi(h^{(m)}(s))
\right] =O(|s|^{m+2}).
\end{displaymath}
\end{prop}
\begin{proof}
The first term on the left-hand side satisfies
\begin{align*}
\overline\partial_1  \omega(s)
 & \cdot 
	( f_j^{(m)}(s) + g_{j\mid (m+1)}(s) + \theta_{j\mid (m+1)}(s) )  \\
 & + R_2 ( \omega(s) \cdot (f_j^{(m)}(s) + g_{j\mid (m+1)}(s) ) ) 
  =O(|s|^{m+2}),
\end{align*}
as we have seen in the proof of Proposition~\ref{prop:10-1}.  
By the construction of $\phi(t)$ (equation~\eqref{eqn:del-B-i-f}), we
have 
\begin{displaymath}
\overline\partial_1\phi(h^{(m)}(s))
 + \overline\partial_1N\overline\partial_1^{\ast} L R_2(\phi(h^{(m)}(s)))=0.
\end{displaymath}
Taking the difference of the last two equations implies
\begin{align*}
\overline\partial_1 & 
	\left[
	\omega(s) \cdot (f_j^{(m)}(s) + g_{j\mid (m+1)}(s)
				+\theta_{j\mid (m+1)}(s) )
	- \phi(h^{(m)}(s))
	\right] \\
& + \overline\partial_1 N \overline\partial_1^{\ast} L
  \left[
	R_2(\omega(s) \cdot 
	  ( f_j^{(m)}(s) + g_{j\mid (m+1)}(s) + \theta_{j\mid (m+1)}(s) ) )\right.\\
& \left.
	\quad - R_2(\phi(h^{(m)}(s)))
  \right]  \\
& + (1-\overline\partial_1N\overline\partial_1^{\ast} L)
  ( R_2(\omega(s) \cdot (f_j^{(m)}(s) + g_{j\mid (m+1)}(s)
			 + \theta_{j\mid (m+1)}(s) ) )
  ) \\
 & \quad =O(|s|^{m+2}).
\end{align*}
The proposition then follows from
Lemmas~\ref{lemma:10-2}~and~\ref{lemma:10-3}.  
\end{proof}

\begin{proof}[Proof of Theorem~\ref{thm:10-2}]
We wish to solve equation~\eqref{eqn:10-21}, which can be written as 
\begin{equation}
\label{eqn:10-29}
\begin{aligned}
\omega(s)  \cdot & ( f_j^{(m)}(s) + g_{j\mid (m+1)}(s) ) 
	- \phi(h^{(m)}(s)) \\
   & = - \overline\partial_{T'} \zeta_{(m+1)}(s) + h_{(m+1)}(s) 
     +O(|s|^{m+2}).
\end{aligned}
\end{equation}
We begin by solving 
\begin{equation}
\label{eqn:10-29a}
\begin{aligned}
\omega(s) \cdot & ( f_j^{(m)}(s) + g_{j\mid (m+1)}(s) 
		  + \theta_{j\mid (m+1)}(s)) 
	- \phi(h^{(m)}(s))  \\
 & = - \overline\partial_{T'} \eta_{(m+1)}(s) + h_{(m+1)} (s)
   +O(|s|^{m+2})
\end{aligned}
\end{equation}
for $\eta_{(m+1)}$ and $h_{(m+1)}$.  
By Proposition~\ref{prop:10-2}, the left-hand side of this equation is
in the kernel of $\overline \partial_1$, modulo $O(|s|^{m+2})$.  The
decomposition of Theorem~\ref{thm:Hodge-Kodaira} implies that 
$\overline \partial_{T'} \eta_{(m+1)}$ and $h_{(m+1)}$, defined as follows,
satisfy equation~\eqref{eqn:10-29a}: 
\begin{align*}
\overline \partial_{T'} \eta_{(m+1)}(s)
       &= - \kappa_{m+1}\left[ 
	D D^* N
	 \left( \omega(s)\cdot (f_j^{(m)}(s)+g_{j\mid (m+1)}(s)\right.\right.\\
&\left.\left.
	    		  \quad + \theta_{j\mid (m+1)}(s) ) 
	   - \phi(h^{(m)}(s)) 
	 \right)
	\right] \\
h_{(m+1)} 
&= \kappa_{m+1} \left[
	H \left( \omega(s)\cdot (f_j^{(m)}(s) + g_{j\mid (m+1)}(s)\right.\right.\\
&\left.\left.\quad
 			    + \theta_{j\mid (m+1)}(s) )
	    - \phi(h^{(m)}(s)) 
	  \right)
	\right].
\end{align*}
Since $DD^* = \overline \partial_0 \rho \rho^* \overline\partial_0^*$,
and $\overline \partial_{T'} = \overline\partial_0$ for elements of $H_0
\subset \Gamma(M,T')$, we may define $\eta_{(m+1)}$ locally by 
\begin{align*}
\eta_{j\mid (m+1)}(s) 
&= - \kappa_{m+1} \left[
	\rho \rho^* \overline\partial_0^* N
	 \left( \omega(s)\cdot (f_j^{(m)}(s)+g_{j\mid (m+1)}(s)\right.\right.\\
&\left.\left.
  \qquad + \theta_{j\mid (m+1)}(s) ) 
	   - \phi(h^{(m)}(s)) 
	\right) 
	\right].  
\end{align*}
To solve equation~\eqref{eqn:10-29}, and thus
equation~\eqref{eqn:10-21}, we simply set 
$\zeta_{(m+1)}=\theta_{(m+1)}+\eta_{(m+1)}$.
This $\zeta$ and $h_{(m+1)}$ solve equation~\eqref{eqn:10-21}, so we
have proved Theorem ~\ref{thm:10-2}.
\end{proof}

Continuing the proof of Theorem~\ref{thm:10-1}, we turn to the proof
of convergence of the formal series.  This part of the proof uses the
standard method of Kodaira and Spencer (see 
\cite{Akahori:1978-family,Akahori-Miyajima:1980}).  We define a Sobolev
$(0,l)$ norm on a power series by setting 
\begin{displaymath}
\norm{ f_j(s) }_{0,l} 
 	= \sum_{|\alpha|=0}^{\infty} 
	  \norm{ f_{j\mid \alpha} }_{0,l} s^\alpha
\qquad \text{ and } \qquad 
\norm{ h(s) }_{0,l}
	= \sum_{|\alpha|=0}^{\infty} 
	  \norm{ h_{\mid \alpha} }_{0,l} s^\alpha.
\end{displaymath}
Consider the power series
\begin{equation}
\label{eqn:A(s)-defn}
A(s) = \frac{b}{16c} 
		\sum_{|\alpha| = 1}^{\infty}
		\left( c^{|\alpha|} / |\alpha|^2 \right) s^{\alpha};
\end{equation}
this converges for any positive $c$.  Moreover, for positive $b$, we
have $A(s)^2 \ll (b/c) A(s)$, where $\ll$ means every coefficient of
the left-hand side is less than the corresponding coefficient of the
right-hand side.  This implies that $A(s)^k \ll (b/c)^{k-1} A(s)$ for
all integers $k \ge 2$.  By choosing suitable $b$ and $c$ (see 
\cite[pp.\ 842--846]{Akahori:1978-family} or 
\cite[Section 3(II), p.\ 832]{Akahori-Miyajima:1980}), 
we wish to show, for any integer $l \ge 3$, that
\begin{equation}
\label{eqn:conv-claim}
\norm{ f_j(s) - z_j(s)}_{0,l} \ll A(s) 
\qquad \text{ and } \qquad 
\norm{ h(s) }_{0,l} \ll A(s).
\end{equation}
(The reason for subtracting $z_j(s)$ in \eqref{eqn:conv-claim}
is because $A(s)$ has no $s^0$ term.)
By the Sobolev embedding theorem, this would give us all the
convergence and regularity claimed in Theorem~\ref{thm:10-1}.

Proof of the convergence~\eqref{eqn:conv-claim} is done by induction
on the partial sums.  That is, we assume that we have 
\begin{equation}
\label{eqn:conv-claim-partial}
\norm{ f^{(m)}_j(s)- z_j(s) }_{0,l} \ll A(s) 
\qquad \text{ and } \qquad 
\norm{ h^{(m)}(s) }_{0,l} \ll A(s),
\end{equation}
then establish the same inequality for $m+1$.  The special properties
of $A(s)$ are used here:  we bound the $m+1$st degree terms with lower
degree terms that we have previously bounded.  As $b$ and $c$ are
chosen properly, we can bound sums of powers of $A(s)$ by $A(s)$
itself.   

The $h_{(m+1)}$ term is well-behaved: for any $l$ there is a constant 
$C_{l}$ such that the harmonic projector $H$ satisfies the estimate 
\begin{displaymath}
\norm{ Hf }_{0,l} \leq C_{l} \norm{ f }
\end{displaymath}
for any $f\in\Gamma(M,E_1)$.  However, we may have to correct $\zeta$
to ensure convergence, because our construction of 
of $\theta_{(m+1)}$ involved first
derivatives of $f_j^{(m)}(s)$.  
Recall from Theorem~\ref{thm:10-2} that
$\zeta_{(m+1)}$ is a solution to 
equation \eqref{eqn:10-21}, which can be viewed as a linear
$\overline\partial_{T'}$ equation for the standard
deformation complex \eqref{eqn:standard-complex}.  
Because $T'$ is a holomorphic vector bundle,
by the results of
\cite{Tanaka:book} there is a Neumann operator
$N_{T'}\colon \Gamma_2(M,T' \otimes ({}^0T'')^*)
\to \Gamma_2(M,T' \otimes ({}^0T'')^*)$ satisfying 
$u = H_{T'}u + \square_{T'} N_{T'}u$ for all
$u$, 
where $\square_{T'} = \overline\partial_{T'}\overline\partial_{T'}^{\ast}+
\overline\partial_{T'}^{\ast}\overline\partial_{T'}$ and
$H_{T'}$ is the projection onto $\ker \square_{T'}$.
Arguing as in 
\cite{Akahori:1978-family}, we let 
\begin{displaymath}
\zeta = - \kappa_{m+1} \left[ 
	 \overline\partial_{T'}^{\ast} N_{T'}
		\left(
    		\omega(s) \cdot (f_j^{(m)}(s) 
				   + g_{j\mid (m+1)} 
	    			   + \theta_{j\mid (m+1)} )
	 	- \phi(h^{(m)}(s) ) 
		\right)
	\right]
\end{displaymath}
and 
\begin{displaymath}
f_j^{(m+1)}(s)
 = f_j^{(m)}(s) + g_{j\mid (m+1)} +\zeta_{j \mid (m+1)}.
\end{displaymath}
It is true that there is a first derivative of $f_j^{(m)}$
in
$\omega(s)\cdot (f_j^{(m)}(s) + g_{j\mid (m+1)} + \theta_{j\mid (m+1)} )
 - \phi(h^{(m)}(s))$, 
but only in the $\C\otimes H$ direction.  
In fact, we recall that 
$\omega(s)\cdot f_j^{(m)}(s)$ is defined on $U_j$ 
by
\begin{displaymath}
(\overline X+\omega(s)\cdot f_j^{(m)}(s)(\overline X))f_j^{(m)}(s) = 0, 
\ \overline X \in {}^0T^{\prime\prime}.
\end{displaymath}
(The CR structure defined on $U_j$ by $\omega(s)\cdot f_j^{(m)}(s)$ 
makes sense globally, modulo $O(|s|^{m+1})$.)  Thus
\begin{align*}
\omega(s)\cdot 
 & f_j^{(m)}(s)z_j(s) + \omega(s)\cdot f_j^{(m)}(s)
	(f_j^{(m)}(s)-z_j(s)) \\
 &+ {\overline X} f_j^{(m)}(s) =O(|s|^{m+2}).
\end{align*}
By the inductive hypothesis, we have
\begin{align*}
\omega(s)\cdot 
 & f_j^{(m)}(s)z_j(s) + \phi(h^{(m)}(s))
	(f_j^{(m)}(s)-z_j(s)) \\
 & + {\overline X} f_j^{(m)}(s) =O(|s|^{m+2}),
\end{align*}
while $\phi(h^{(m)}(s))$ takes its values in $^0T^{\prime}$.  
Since the composition $\overline\partial_{T'}^{\ast}N_{T'}$ 
of the adjoint operator and the standard Neumann
operator gains 1 in 
this direction, there 
is no problem in the convergence of our formal solution.
This finishes the proof of Theorem~\ref{thm:10-1}.
\end{proof}

\end{document}